%% file: betaMixingEVT.tex
\documentclass[a4paper,12pt]{article}
\renewcommand{\baselinestretch}{1.5}
\usepackage{amsmath}
\usepackage{amssymb}
\usepackage{color}

\usepackage{epsfig}
\usepackage{lscape}
\usepackage{enumitem}
\usepackage{float}

\usepackage[english]{babel}
\usepackage{booktabs}
\usepackage{bm}
\usepackage{bbm}
\usepackage{cases}
\usepackage{eucal}
\usepackage{float}
\usepackage{dsfont}
\usepackage{epsfig}
\usepackage{fontenc}

\usepackage{mathabx}

\def\ind{ {{\rm 1}\hskip-2.2pt{\rm l}}}
 \usepackage{url}
 \usepackage[hmargin = 1.6cm, vmargin = 3cm]{geometry}
\usepackage[colorlinks, allcolors = blue]{hyperref}

\usepackage{lscape}
\usepackage{mathrsfs}
\usepackage{natbib}
\usepackage{rotating}
\usepackage{subfig}
\usepackage{verbatim}
\usepackage[table]{xcolor}

\input{colors.tex}
\input{environments.tex}
\input{header.tex}

\input{tools.tex}

\renewcommand{\baselinestretch}{1.5}

\newcommand{\blind}{1}
  \date{}

\begin{document}

\bibliographystyle{agsm}  
\def\spacingset#1{\renewcommand{\baselinestretch}%
{#1}\small\normalsize} \spacingset{1}

\if1\blind
{
  \title{Risk measure estimation for $\beta$-mixing time series and applications}
\author{{Val\'erie \textsc{Chavez-Demoulin}, Armelle \textsc{Guillou}}
\thanks{Val\'erie Chavez-Demoulin is Professor, Faculty of Business and Economics (HEC), Universit\'e de Lausanne, Switzerland (\href{mailto:valerie.chavez@unil.ch}{\it valerie.chavez@unil.ch}). Armelle Guillou is Professor, Universit\'e de Strasbourg \& CNRS, IRMA UMR 7501, France (\href{mailto:armelle.guillou@math.unistra.fr}{\it armelle.guillou@math.unistra.fr}). The research was partially funded by a research grant (VKR023480) from VILLUM FONDEN. We thank Laurens de Haan, C\'ecile Mercadier and Chen Zhou for several discussions about the proof of Theorem 4.2 in \cite{dHMZ2016} which has allowed to detect an error in their paper. We thank them for having provided the authorization to mention the corrected estimator (\ref{estimcorrige}).}}
  \maketitle
} \fi

\if0\blind
{
  \bigskip
  \bigskip
  \bigskip
  \begin{center}
    {\LARGE\bf Risk measure estimation for $\beta$-mixing time series and applications}
\end{center}
  \medskip
} \fi
\bigskip
\begin{abstract}
In this paper, we discuss the application of extreme value theory in the context of stationary $\beta$-mixing sequences that belong to the Fr\'echet domain of attraction. In particular, we propose a methodology to construct bias-corrected tail estimators. Our approach is based on the combination of two estimators for the extreme value index to cancel the bias. The resulting estimator is used to estimate an extreme quantile. In a simulation study, we outline the performance of our proposals that we compare to alternative estimators recently introduced in the literature. Also, we compute the asymptotic variance in specific examples when possible. Our methodology is applied to two datasets on finance and environment.
\end{abstract}
\vfill

\newpage
\spacingset{1.4} 

\section{\textsf{INTRODUCTION}}
\label{intro}
Quantitative Risk Management (QRM) has become an inevitable field with aimed at building models to understand the risks of financial portfolios and environmental hazards. For the most complete treatment of the theoretical concepts and modeling tools of QRM, we refer to the book by \cite{mcneil2015}. Building such models is now a crucial task across the banking and insurance industries under the regulatory obligations of the Basel Committee
on Banking Supervision and Solvency 2. In financial risk management, research related to regulatory risk measures such as the Value-at-Risk (VaR) has received lot of attention over the last decades. 
In environmental risk management, and due to the increasing frequency of extreme events \citep{swissre2014, EKR2016} and their disastrous societal impact, estimating risk measures such as the return level is of vital importance. Both these risk measures (VaR and return level) rely on high quantile estimation in the tail region of the observations distribution. In this context of tail modeling, extreme 
value theory (EVT) offers strong and adequate statistical tools. Classical EVT models are based on the independent and identically distributed (i.i.d.) assumption. This assumption is however very often violated in practice. Financial time series, for instance, show volatility clustering and environmental data typically exhibit serial dependence. Many EVT-related papers for time series address the modeling of such features \citep[see, for instance,][]{mcneil2000,CES2014}.
In this paper, we introduce a new asymptotically unbiased high quantile estimator for stationary time series 
such as the very commonly used autoregressive (AR), the moving average (MA) and the generalized autoregressive conditional heteroskedasticity (GARCH) models. More precisely, we propose a new estimator of high quantiles for $\beta$-mixing stationary time series with heavy-tailed distribution. The estimator is based on an asymptotically unbiased estimator of the extreme value index. The advantage of our estimator is twofold: first, it directly handles the serial feature of such $\beta$-mixing time series contrary to other methods that need a pre-filtering of the heteroskedasticity before applying the standard estimator. As an example, we refer to the two step-method of \cite{mcneil2000}. Second, it improves the alternative bias correction procedure proposed by \cite{dHMZ2016} for $\beta$-mixing series. 

Throughout the paper, we assume that $(X_1, X_2, ...)$ is a $\beta$-mixing time series, that is, a series such that
$$\beta(m):=\sup_{p\geq 1} \mathbb E\left\{\sup_{C\in {\cal B}^\infty_{p+m+1}} \left|\mathbb P(C|{\cal B}^p_1)-\mathbb P(C)\right|\right\} \longrightarrow 0,$$
as $m\to \infty$, where ${\cal B}_i^j$ denotes the $\sigma-$algebra generated by $X_i, ..., X_j$. Loosely speaking, $\beta(m)$ measures the total variation distance between the unconditional distribution of the future of the time series and the conditional distribution of the future given the past of the series when both are separated by $m$ time points. 
Let $F$ be the common marginal distribution function of $X_i, i\in \mathbb N$, which is assumed to belong to the Fr\'echet domain of attraction, that is, the tail quantile function $U:=(1/(1-F))^\leftarrow$  where $^\leftarrow$ denotes the left continuous inverse function, satisfies
\begin{eqnarray}
\lim_{t\to \infty} {U(tx)\over U(t)} = x^\gamma, \quad \forall x>0.
\label{tailU}  
\end{eqnarray}
The estimation of the extreme value index $\gamma$ has been extensively studied in the case of i.i.d. random variables, but only few papers consider this topic in case of time series with serial dependence features. We can mention, among others, \cite{hsing1991}, \cite{drees2000} and \cite{drees2003} and very recently \cite{dHMZ2016}. As in the i.i.d. context, the simplest estimator for $\gamma>0$ is the Hill estimator \citep{hill1975} defined as
$$\widehat \gamma^H_k:={1\over k}\sum_{i=1}^k \log X_{n-i+1,n}-\log X_{n-k,n},$$
where $X_{1,n}\leq \cdots \leq X_{n,n}$ denote the order statistics and $k$ is an intermediate sequence, that is, a sequence such that $k\to \infty$ and $k/n\to 0$ as $n\to \infty$.

To prove the asymptotic normality of a tail parameter such as the Hill estimator, we need a second order condition which specifies the rate of convergence for the left-hand side in (\ref{tailU}) to its limit. This condition can be formulated in different ways, below we state it in terms of the logarithm since it is this formulation that we will use later.

{\bf Second order condition $(C_{SO})$.} {\it Suppose that there exists a positive or negative function $A$ with $\lim_{t\to \infty} A(t)=0$ and a real number $\rho<0$ such that
\begin{eqnarray*}
\lim_{t\to \infty} {\log U(tx)- \log U(t)-\gamma \log x \over A(t)} = {x^\rho-1\over \rho}, \quad \forall x>0.
\end{eqnarray*} 
}

The rate of convergence for the function $A$ to 0 is crucial if we want to exhibit the bias term of the estimator of a tail parameter. Under the assumption that the intermediate $k-$sequence is such that $\sqrt k A(n/k)\to \lambda\in \mathbb R$, and assuming the following
regularity conditions on the $\beta$-mixing coefficients :\\
{\bf Regularity conditions $(C_R)$.} {\it There exist $\varepsilon>0$, a function $r$ and a sequence $\ell_n$ such that, as $n \to \infty$,
\begin{enumerate}[label=(\alph*)]
\item ${\beta(\ell_n) \over \ell_n} n + \ell_n {\log^2 k\over \sqrt k} \longrightarrow 0;$
\item ${n \over \ell_n k} \mbox{Cov}\left(\sum_{i=1}^{\ell_n} \ind_{\{X_i>F^\leftarrow(1-kx/n)\}}, \sum_{i=1}^{\ell_n} \ind_{\{X_i>F^\leftarrow(1-ky/n)\}}\right)\longrightarrow r(x, y),\quad \forall\, 0\leq x, y\leq 1+\varepsilon;$
\item For some constant $C$:
$${n \over \ell_n k} \mathbb E\left[\left(\sum_{i=1}^{\ell_n} \ind_{\{F^\leftarrow(1-ky/n)<X_i\leq F^\leftarrow(1-kx/n)\}}\right)^4\right] \leq C(y-x), \quad \forall\,  0\leq x< y\leq 1+\varepsilon \mbox{ and } n\in \mathbb N,$$
\end{enumerate}
}
\cite{drees2000} has established the asymptotic normality of $\widehat \gamma^H_k$
\begin{eqnarray}
\sqrt k(\widehat \gamma^H_k-\gamma) \stackrel{d}{\longrightarrow} {\cal N}\left({\lambda \over 1-\rho}, \sigma^2\right),
\label{Hill}
\end{eqnarray}
where $\sigma^2$ depends on the covariance structure $r$, but has a simple expression in the i.i.d. context, where it is equal to $\gamma^2$.
In practice, the bias term of $\widehat \gamma^H_k$ can be important depending on whether $\rho$ is close to zero or not, since under the second order condition $(C_{SO})$, the function $|A|$ is regularly varying at infinity with index $\rho$. This explains all the literature spread on bias correction in the i.i.d. context. On the contrary, in case of stationary $\beta$-mixing time series only the very recent paper by \cite{dHMZ2016} deals with this problem and proposes a bias-corrected estimator for $\gamma$. Their method consists first in estimating the bias term of $\widehat \gamma^H_k$ and second in subtracting it from $\widehat \gamma^H_k$. A similar approach is also used in their paper to estimate a high quantile $x_p=U(1/p)$ with $p\to 0$.

The procedure we propose in this paper is an alternative approach to construct bias-corrected tail estimators. First, we introduce a class of  estimators for $\gamma$ which can be viewed as statistical tail functionals, $T(Q_n)$, where $Q_n$ is the tail quantile function defined as $Q_n(t):=X_{n-[kt],n}, 0<t<n/k$, and $T$ is a suitable functional. Then, we combine two estimators for $\gamma$ of this class to cancel the asymptotic bias term. The resulting unbiased estimator for $\gamma$ can then be used to construct an asymptotically unbiased estimator of a high quantile. 

The paper is organized as follows: our approach is described in details in Section~\ref{construction}. Section~\ref{simulation} presents, using some examples, the finite sample performance of our high quantile estimator based on simulation studies. Two real data applications illustrate the use of our estimator in Section~\ref{applications}: one in the financial context of market risk data and the other in the environmental situation of hourly wind speed data. We conclude in Section \ref{conclusion}. All the related theoretical proofs are detailed in the appendix. 

\section{\textsf{DESCRIPTION OF OUR METHODOLOGY}}
\label{construction}
\cite{GG2013} introduced a class of weighted function estimators for the tail dependence coefficient $\eta$ in the bivariate extreme value framework. Combining two of their estimators, they are able to construct an asymptotically unbiased estimator for $\eta$. In this paper, we propose to use the same methodology in the case of $\beta$-mixing sequences to estimate a tail parameter such as the extreme value index or an extreme quantile.

\subsection{\textsf{Estimation of the extreme value index}}
\label{sectiongamma}
For any measurable function $z: [0, 1] \to \mathbb R$, we consider the functional
\begin{eqnarray*}
T_K(z)=\left\{
    \begin{array}{ll}
   \displaystyle{ \int_0^1 \log {z(t)\over z(1)} d(tK(t)) }& \mbox{if the right-hand side is defined and finite,}\\
    0 & \mbox{otherwise.}
     \end{array}
\right.
\end{eqnarray*}
This leads to the following class of estimators for $\gamma$:
\begin{eqnarray*}
\widehat \gamma_k(K)=T_K(Q_n)=\int_0^1 \log {Q_n(t)\over Q_n(1)} d(tK(t)),
\end{eqnarray*}
where $K$ is a function with support on $(0, 1)$. Some assumptions on $K$ are required if we want to derive the asymptotic normality of our class of estimators. They can be formulated as:\\
{\bf Assumption $(C_K)$.} {\it Let $K$ be a function such that $\int_0^1 K(t)dt=1$. Suppose that $K$ is continuously differentiable on $(0, 1)$ and that there exist $M>0$ and $\tau \in [0, 1/2)$ such that $|K(t)|\leq M \, t^{-\tau}$.}

These conditions are not restrictive but are satisfied by the usual weight functions used in the literature, including the power kernel $K(u)=(1+\nu)u^\nu, \nu\geq 0,$ and the log-weight function $K(u)=(-\log u)^\nu/\Gamma(1+\nu), \nu\geq 0$. In particular, we note that the classical Hill estimator $\widehat \gamma^H_k$ can be viewed as a particular case of our power kernel-type estimator corresponding to $\nu=0$: $\widehat \gamma^H_k=\widehat \gamma_k(K)$ with $K(u)=1.$

The aim of Theorem~1 is to provide the asymptotic normality of our class of estimators with the explicit bias term.
\begin{theorem} Let $(X_1, X_2, ...)$ be a stationary $\beta$-mixing time series with a continuous common marginal distribution function $F$ and assume $(C_{S0})$, $(C_R)$ and $(C_K)$. Suppose that $k$ is an intermediate sequence such that $\sqrt k A(n/k)\to \lambda\in \mathbb R$. We have
\begin{eqnarray*}
\sqrt k \left\{\widehat \gamma_k(K)-\gamma-A\left({n\over k}\right)\int_0^1t^{-\rho}K(t)dt\right\} \stackrel{d}\longrightarrow \gamma\int_0^1 \left[t^{-1}W(t)-W(1)\right] d(tK(t))
\end{eqnarray*}
where $(W(t))_{t\in [0, 1]}$ is a centered Gaussian process with covariance function $r$ defined in $(C_R)$.
\label{thm1}
\end{theorem}
This result is similar to Theorem 1 in \cite{GG2013} in the i.i.d. context, although in the latter the centered Gaussian process $(W(t))_{t\in [0, 1]}$ is a standard Brownian motion. \\

From the covariance structure, we deduce in the next corollary the asymptotic normality of $\widehat \gamma_k(K)$ and its asymptotic variance.

{\bf Corollary 1.} {\it Under the assumption of Theorem~1, we have
\begin{eqnarray*}
\sqrt k \{\widehat \gamma_k(K)-\gamma\} \stackrel{d}{\longrightarrow}
{\cal N}\left(\lambda{\cal AB}(K), {\cal AV}(K)\right)
\end{eqnarray*}
where 
\begin{eqnarray*}
{\cal AB}(K)&=&\int_0^1 t^{-\rho}K(t)dt\\
\mbox{and }{\cal AV}(K)&=&\gamma^2\int_0^1\int_0^1 \left[{r(t, s)\over ts} - {r(t, 1)\over t} - {r(1, s) \over s} +r(1, 1)\right] d(tK(t)) d(sK(s)).
\end{eqnarray*}
}

{\bf Some remarks.}
\begin{itemize}
\item If we assume that the $X_i$s are i.i.d. and that $K(t)=1$ for $t\in (0, 1)$, then Corollary 1 gives the asymptotic normality of the Hill estimator (\ref{Hill}) with $\sigma^2=\gamma^2$.
\item If we assume that the $X_i$s are i.i.d. but nothing on the function $K$ (except that it satisfies condition $(C_K)$), then the asymptotic variance given in Corollary 1 can be reduced to
$${\cal AV}(K)=\gamma^2\left\{\int_0^1\int_0^1 {\min(s,t)\over ts} d(tK(t)) d(sK(s))-K^2(1)\right\}.$$
\end{itemize}

Now, our aim is to propose an asymptotically unbiased estimator for $\gamma$. For this aim, we propose to use two functions $K_1$ and $K_2$ satisfying $(C_K)$ and to consider a
mixture of them in the form $K_\Delta(t)=\Delta K_1(t)+(1-\Delta)K_2(t)$ for $\Delta \in \mathbb R$. Clearly $K_\Delta$ also satisfies condition $(C_K)$ and hence by Corollary 1, the asymptotic bias of this new estimator $\widehat \gamma_k({K_\Delta})$ is given by
$${\lambda \over \sqrt k} \, {\cal AB}(K_\Delta)={\lambda \over \sqrt k} \, \int_0^1 t^{-\rho} K_\Delta(t)dt={\lambda \over \sqrt k} \, \left\{\Delta {\cal AB}(K_1)+(1-\Delta) {\cal AB}(K_2)\right\}.$$
Equating the right-hand side of the above equation to zero leads to the value of $\Delta$ eliminating the asymptotic bias
\begin{eqnarray}
\Delta^*={{\cal AB}(K_2) \over {\cal AB}(K_2)-{\cal AB}(K_1)} \quad \mbox{provided ${\cal AB}(K_1)\not = {\cal AB}(K_2)$.}
\label{deltaopt}
\end{eqnarray}
This result is formalized in the next corollary where $\widehat \gamma_k(K_{\Delta^*})$ is shown to be asymptotically unbiased in the sense that the mean of its limiting distribution is zero, whatever the value of $\lambda$.

{\bf Corollary 2.} {\it Under the assumptions of Theorem~1 and assuming that $K_1$ and $K_2$ satisfy condition $(C_K)$ with ${\cal AB}(K_1)\not = {\cal AB}(K_2)$, we have
\begin{eqnarray*}
\sqrt k \{\widehat \gamma_k(K_{\Delta^*})-\gamma\} \stackrel{d}{\longrightarrow}
{\cal N}\left(0, {\cal AV}(K_{\Delta^*})\right).
\end{eqnarray*}
}
An open problem is to determine whether among this class of unbiased estimators we can find the asymptotically unbiased estimator with minimum variance.
This question is solved in the i.i.d. framework under a slightly stronger condition than $(C_K)$ \citep[see Theorem 2 and Corollary 4 in][]{GG2013}
where the ``optimal'' function is given by
 \begin{eqnarray}
K_{\Delta^*_{opt}}(t)=\left({1-\rho\over \rho}\right)^2 - {(1-\rho)(1-2\rho)\over \rho^2} t^{-\rho}, \quad t\in (0, 1).
\label{Kopt}
\end{eqnarray}
Note that this function can be viewed as a mixture between  $K_1(t):=1$ and $K_{2,\rho}(t):=(1-\rho)t^{-\rho}$ with $t\in (0, 1)$ and $\Delta^*$ as in (\ref{deltaopt}).  In that case,
the minimal variance is given by
\begin{eqnarray}
{\cal AV}(K_{\Delta^*_{opt}})=\gamma^2\left({1-\rho\over \rho}\right)^2.
\label{Avaropt}
\end{eqnarray}
Although we are not able to show that this interesting property is preserved in this new context, we recommend using this ``optimal'' function also in case of $\beta$-mixing sequences. However, from a practical point of view, this
cannot be done directly since $\rho$ is unknown. To solve this issue, two natural options can be proposed, either to replace $\rho$ by a canonical choice, or by an external estimator.

The aim of the next corollary is to give the asymptotic normality of our class of estimators for $\gamma$ in case $\rho$ is replaced by some fixed value $\widetilde \rho$.

{\bf Corollary 3.} {\it  Let $(X_1, X_2, ...)$ be a stationary $\beta$-mixing time series with a continuous common marginal distribution function $F$ and assume $(C_{S0})$ and $(C_R)$. Suppose that $k$ is an intermediate sequence such that $\sqrt k A(n/k)\to \lambda\in \mathbb R$. We have
 \begin{eqnarray*}
\sqrt k \{\widehat \gamma_k(K_{\widetilde\Delta^*_{opt}})-\gamma\} \stackrel{d}{\longrightarrow}
{\cal N}\left(\lambda {(1-\widetilde \rho)(\widetilde \rho-\rho) \over \widetilde \rho(1-\rho)(1-\widetilde\rho-\rho)}, {\cal AV}(K_{\widetilde \Delta^*_{opt}})\right),
\end{eqnarray*}
where $K_{\widetilde\Delta^*_{opt}}$ is defined as $K_{\Delta^*_{opt}}$ in (\ref{Kopt}) with $\rho$ replaced by $\widetilde \rho$. 
}

Although one clearly loses the bias correction, the extreme value index estimators are not very sensitive to such a misspecification and thus our estimator $\widehat \gamma_k(K_{\widetilde\Delta^*_{opt}})$ can still outperform the estimators that are not corrected for bias. Note also that, as expected, if $\widetilde\rho=\rho$, we recover Corollary 2. However, to keep the asymptotically unbiased property, we can also replace $\rho$ by an external estimator $\widehat \rho_{k_\rho}$, consistent in probability, which depends on an intermediate sequence $k_\rho$. This leads to the following general result.
\begin{theorem} Let $(X_1, X_2, ...)$ be a stationary $\beta$-mixing time series with a continuous common marginal distribution function $F$ and assume $(C_{S0})$ and $(C_R)$. Let $\widehat \rho_{k_\rho}$ be an external estimator for $\rho$, consistent in probability, which depends on an intermediate sequence $k_\rho$. If $k$ is an intermediate sequence such that $\sqrt k A(n/k)\to \lambda\in \mathbb R$, then we have
\begin{eqnarray*}
\sqrt k \{\widehat \gamma_k(K_{\widehat\Delta^*_{opt}})-\gamma\} \stackrel{d}{\longrightarrow}
{\cal N}\left(0, {\cal AV}(K_{\Delta^*_{opt}})\right),
\end{eqnarray*}
where $K_{\widehat\Delta^*_{opt}}$ is defined as $K_{\Delta^*_{opt}}$ in (\ref{Kopt}) with $\rho$ replaced by $\widehat \rho_{k_\rho}$. 
\label{thm2}
\end{theorem}

Note that this theorem cannot be viewed as a consequence of one of the two first corollaries because $K_{\widehat \Delta^*_{opt}}$ depends on $\widehat \rho_{k_\rho}$ and thus on $n$ which means that Corollary 1 cannot be used directly. Moreover, it cannot be written as a mixture of the form $\Delta^*K_1+(1-\Delta^*)K_{2,\widehat \rho_{k_\rho}}$ with $\Delta^*$ defined by (\ref{deltaopt}) since this expression would depend on both $\widehat \rho_{k_\rho}$ and $\rho$, which is unknown. This implies that this theorem cannot be viewed as the counterpart of Proposition 2 in \cite{GG2013} where the limiting distribution of the normalized bias-corrected estimator of the tail dependence coefficient is established, but only in cases where the two kernels are assumed to be independent on $\rho$. Such a result can also be obtained in our framework, but we omit it since  we recommend here to use our ``optimal" function $K_{\widehat\Delta^*_{opt}}$. Note that \cite{GG2008} mention a result similar to our Theorem~2 in their framework as an open problem.

A possible choice for $\widehat \rho_{k_\rho}$ is that proposed by \cite{GdHP2002}, and also used in \cite{dHMZ2016}:
\begin{eqnarray}
\widehat \rho_{k}:={-4+6\, S_{k}^{(2)}+\sqrt{3\, S_{k}^{(2)}-2} \over 4\, S_{k}^{(2)} -3} \quad\mbox{provided $S_{k}^{(2)}\in \left({2 \over 3}, {3 \over 4}\right)$},
\label{rhoestimateur}
\end{eqnarray}
where
\begin{eqnarray*}
S_{k}^{(2)}:={3 \over 4}\, {\left[M_{k}^{(4)}-24\left(M_{k}^{(1)}\right)^4\right]\left[M_{k}^{(2)}-2\left(M_{k}^{(1)}\right)^2\right] \over \left[M_{k}^{(3)}-6\left(M_{k}^{(1)}\right)^3\right]^2}
\end{eqnarray*}
with
$$M_{k}^{(\alpha)}:= {1 \over k} \sum_{i=1}^{k} \left(\log X_{n-i+1,n}-\log X_{n-k,n}\right)^{\alpha}, \alpha\in \mathbb N. $$

In that case, we have the following corollary.

{\bf Corollary 4.} {\it Let $(X_1, X_2, ...)$ be a stationary $\beta$-mixing time series with a continuous common marginal distribution function $F$ and assume $(C_{S0})$ and $(C_R)$. Let $\widehat \rho_{k_\rho}$ be the external estimator for $\rho$ defined in (\ref{rhoestimateur}) where the intermediate sequence $k_\rho$ satisfies $\sqrt {k_\rho} A\left({n \over k_\rho}\right) \to \infty$.
If $k$ is another intermediate sequence such that $\sqrt k A(n/k)\to \lambda\in \mathbb R$, then we have
\begin{eqnarray*}
\sqrt k \{\widehat \gamma_k(K_{\widehat\Delta^*_{opt}})-\gamma\} \stackrel{d}{\longrightarrow}
{\cal N}\left(0, {\cal AV}(K_{\Delta^*_{opt}})\right).
\end{eqnarray*}
}

Compared to Theorem 4.1 in \cite{dHMZ2016}, the assumptions of our Corollary 4 are less constraining, in particular we do not need a third order condition and only $\sqrt{k_\rho} A\left({n \over k_\rho}\right) \to \infty$ is required on the intermediate sequence $k_\rho$. This is due to the fact that we need only the consistency in probability for the $\rho-$estimator, not its asymptotic normality. However, our rate of convergence $\sqrt k$ is smaller than that obtained by   \cite{dHMZ2016} since they assume $\sqrt k A\left({n \over k}\right)\to \infty$ in place of our condition $\sqrt k A\left({n \over k}\right)\to \lambda$. This implies that, although their rate of convergence also has the form $\sqrt k$, their intermediate sequence is larger than ours, taking into account that the function $|A|$ is regularly varying at infinity with index $\rho$. Despite our slower rate,  our estimator outperforms the one proposed by \cite{dHMZ2016} in finite samples situation.

\subsection{\textsf{Estimation of an extreme quantile}}
The estimation of an extreme value index is in general only an intermediate goal. In practice, we are much more interested in the estimation of an extreme quantile 
\begin{equation}
x_p=U(1/p),
\label{quantile}
\end{equation}
$p\to 0$. As mentioned in Section \ref{intro}, the VaR and the return level are extreme quantiles consisting of standard risk measures of finance and environment, respectively.

In this section, we illustrate the applicability of our methodology in the case of the estimation of an extreme quantile $x_p$.
To understand heuristically the construction of our estimator, we start with our second-order condition $(C_{SO})$, according to which
$${U(tx) \over U(t)} \simeq x^\gamma \exp\left\{A(t){x^\rho-1\over \rho}\right\}.$$
By setting $tx=1/p$ and $t=Y_{n-k,n}$ where $Y_i$ is a random variable from a standard Pareto distribution, since $X_{n-k,n}=U(Y_{n-k,n})$, 
we obtain the following approximation
\begin{eqnarray*}
x_p  &\simeq& X_{n-k,n} \left({1\over pY_{n-k,n}}\right)^\gamma \exp\left\{A\left(Y_{n-k,n}\right) {\left({1\over pY_{n-k,n}}\right)^\rho-1\over \rho}\right\}\\
&\simeq&X_{n-k,n} \left({k \over np}\right)^\gamma \exp\left\{A\left({n \over k}\right){({k \over np})^\rho - 1\over \rho}\right\},
\end{eqnarray*}
where the last step follows from replacing $Y_{n-k,n}$ by its expected value $n/k$. Note that the first part on the right-hand side (except the exponential term)
is exactly a Weissman-type estimator when $\gamma$ is replaced by an estimator \citep[see][]{weiss1978}. Thus this exponential term can be viewed as a correcting term since $A(n/k)$ tends to 0 (and thus the exponential to 1). To take this correcting factor into account, we need to estimate $A(n/k)$. For this aim, we use our Proposition~1 in the appendix according to which, for any $\varepsilon>0$,
\begin{eqnarray*}
\log {Q_n(t)\over Q_n(1)} &=& \log {Q_n(t)\over U({n\over k})} - \log {Q_n(1)\over U({n\over k})}\\
&=& - \gamma \log t +{\gamma \over \sqrt k} \left[t^{-1} W(t)-W(1)\right]+\widetilde A\left({n \over k}\right) {t^{-\rho}-1\over \rho} + {o\left(t^{-{1\over 2}-\varepsilon}\right)\over \sqrt k},
\end{eqnarray*}
where $\widetilde A\sim A$. This implies that
\begin{eqnarray*}
\sqrt k \left\{\widehat \gamma_k(K_1)-\widehat \gamma_k(K_{2,\rho})+A\left({n \over k}\right) {\rho^2\over (1-\rho)(1-2\rho)}\right\} \stackrel{d}{\longrightarrow} \gamma \int_0^1 \left[t^{-1} W(t)-W(1)\right] d\left(t\left[1-K_{2,\rho}(t)\right]\right)
\end{eqnarray*}
which is asymptotically normal ${\cal N}(0, {\cal AV}(K_1-K_{2,\rho}))$. Thus we can approximate
$$A\left({n \over k}\right) {\rho^2\over (1-\rho)(1-2\rho)} \qquad \mbox{by  } -\left[\widehat \gamma_k(K_1)-\widehat \gamma_k(K_{2,\rho})\right]$$
which means that $A\left({n \over k}\right)$ can be estimated by
$$-{(1-\xi)(1-2\xi) \over \xi^2} \left[\widehat \gamma_k(K_1)-\widehat \gamma_k(K_{2,\xi})\right]$$
where $\xi$ can be either a consistent estimator for $\rho$ or a canonical negative value.\\
Our final extreme quantile estimator is then
\begin{eqnarray}
\widehat x_{p,\xi} = X_{n-k,n} \left({k \over np}\right)^{\widehat \gamma_k(K_{\widehat\Delta^*_{opt}})} \exp\left\{-{(1-\xi)(1-2\xi)\over \xi^2} \left[\widehat \gamma_k(K_1)-\widehat \gamma_k(K_{2,\xi})\right] {\left({k\over np}\right)^{\xi}-1\over \xi}\right\}.
\label{quantileestimator}
\end{eqnarray}
The aim of the next theorem is to prove that, under suitable assumptions, this estimator is asymptotically unbiased.
\begin{theorem} Let $(X_1, X_2, ...)$ be a stationary $\beta$-mixing time series with a continuous common marginal distribution function $F$ and assume $(C_{S0})$ and $(C_R)$. Let $\widehat \rho_{k_\rho}$ be an external estimator for $\rho$, consistent in probability, which depends on an intermediate sequence $k_\rho$. Consider now an intermediate sequence $k$ such that $\sqrt k A(n/k)\to \lambda\in \mathbb R$ and assume that
$p=p_n$ such that ${k \over np}\to \infty$, ${\log(np)\over \sqrt k} \to 0$ and $n^{-a} \log p \to 0$ for all $a>0$. Then, we have
\begin{eqnarray*}
{\sqrt k \over \log {k\over np}} \left( {\widehat x_{p, \xi} \over x_p} -1\right)\stackrel{d}{\longrightarrow} {\cal N}\left(0, {\cal AV}(K_{\Delta^*_{opt}})\right), 
\end{eqnarray*}
where $\xi$ is either a canonical negative value $\widetilde \rho$ or an estimator $\widehat \rho$ consistent in probability such that $|\widehat \rho-\rho|=O_\mathbb P(n^{-\varepsilon})$ for some $\varepsilon>0$.
\label{thm3}
\end{theorem}

Note that the assumption $\log(np)/\sqrt k \to 0$ is useful in order to ensure that the rate of convergence for our extreme quantile estimator tends to infinity.
The other condition on $p$, that is, $n^{-a} \log p \to 0$ for all $a>0$, is only technical and not binding. This is usual in the context of extreme quantile estimation, \citep[see, for instance,][]{Matthys2004}. In the  latter, an extreme quantile estimator in the context of i.i.d. censored observations has been proposed, and its asymptotic normality also requires the condition $|\widehat \rho-\rho|=O_\mathbb P(n^{-\varepsilon})$ for some $\varepsilon>0$.

From the proof of Theorem~3, it becomes clear that the exponential term in (\ref{quantileestimator}) does not influence the limiting distribution. However, in finite samples situations, this factor typically leads to improved overall stability of the quantile estimates as a function of $k$.

\section{\textsf{EXAMPLES AND SIMULATIONS}}
\label{simulation}
Our aim in this section is to compare our estimator with the one proposed by \cite{dHMZ2016}.
Unfortunately, without specifying the covariance structure $r$ in $(C_R)$, it is impossible to compare our asymptotic variance ${\cal AV}(K_{\Delta^*_{opt}})$ with that of the asymptotic unbiased estimator of \cite{dHMZ2016} in its full generality. However, in the specific case of i.i.d. observations where we know that $r(s, t)=\min(s, t)$, our asymptotic variance given in (\ref{Avaropt}) is clearly smaller than the asymptotic variance of the extreme value index proposed by \cite{dHMZ2016}, which is
$${\gamma^2\over \rho^2} \left(\rho^2+(1-\rho)^2\right).$$

To complete this comparison, we consider below several models commonly found in practice, with an explicit expression of the covariance structure $r$ for two of the models. This allows us to provide an explicit comparison between the two estimators.

\subsection{\textsf{Autoregressive (AR) model}}

Let $\varepsilon_i$ be i.i.d. variables with a positive Lebesgue density $f_\varepsilon$ which is $L_1-$Lipschitz continuous, that is,
$$\int |f_\varepsilon(\varepsilon+u)-f_\varepsilon(\varepsilon)|d\varepsilon = O(u) \qquad \mbox{as } u \searrow 0.$$
Assume that
$$1-F_\varepsilon(\varepsilon)\sim q \varepsilon^{-1/\gamma}\ell(\varepsilon) \qquad \mbox{and} \qquad F_\varepsilon(-\varepsilon)\sim (1-q) \varepsilon^{-1/\gamma}\ell(\varepsilon)\qquad \mbox{as }\varepsilon\to \infty,$$
for some slowly varying function $\ell$ and $q\in (0, 1)$. Consider now the stationary solution of the AR(1) equation
\begin{equation}
X_i=\theta X_{i-1}+\varepsilon_i, \quad \theta\in (0, 1).
\label{AR}
\end{equation}

The regularity conditions $(C_R)$ hold with 
$$r(x, y)=\min(x, y)+\sum_{m=1}^{\infty} \{c_m(x,y)+c_m(y,x)\}$$
where $c_m(x,y)=\min\left(x, y\theta^{m\over \gamma}\right).$

Direct but tedious computations lead to the following asymptotic variance for our estimator
\begin{eqnarray}
{\cal AV}(K_{\Delta^*_{opt}})=\gamma^2\left({1-\rho \over \rho}\right)^2 r(1, 1),
\label{asympvar}
\end{eqnarray}
that is always smaller than that obtained by \cite{dHMZ2016} under the same framework, which is
$$\sigma^2(\theta, \gamma, \rho):= {\gamma^2 \over \rho^2} \left\{\left[(1-\rho)^2+\rho^2\right]r(1, 1)+2\rho(1-\rho){\theta^{1 \over \gamma} \log \theta^{1 \over \gamma} \over (1-\theta^{1\over \gamma})^2}\right\}.$$
Note also that, compared with the i.i.d. case, our asymptotic variance ${\cal AV}(K_{\Delta^*_{opt}})$ is increased by the factor $r(1, 1)>1$, see (\ref{Avaropt}). In addition, if $K(t)=1$ for $t\in (0, 1)$, our estimator $\widehat \gamma_k(K)$ reduces to the classical Hill estimator $\widehat \gamma^H_k$ and according to \cite{drees2000}  \citep[see also][]{starica1999}, under serial dependence, the asymptotic variance of $\widehat \gamma^H_k$ is $\gamma^2 r(1, 1)$. The latter value is smaller than ${\cal AV}(K_{\Delta^*_{opt}})$, but $\widehat \gamma^H_k$ is not asymptotically unbiased.

 \subsection{\textsf{Moving average (MA) model}}
 
Assume that $\varepsilon_i$ satisfies the same assumptions as for the AR(1) model and consider this time the stationary solution of the MA(1) equation
\begin{equation}
X_i=\theta \varepsilon_{i-1}+\varepsilon_i.
\label{MA}
\end{equation}
In that case, the regularity conditions $(C_R)$ are also satisfied with
$$r(x, y)=\min(x, y)+\left(1+\theta^{1 \over \gamma}\right)^{-1}\left\{\min\left(x, y\theta^{1\over \gamma}\right)+\min\left(y, x\theta^{1\over \gamma}\right)\right\}.$$
Again tedious computations show that the same expression for ${\cal AV}(K_{\Delta^*_{opt}})$ as that given in (\ref{asympvar}) is valid for the MA(1) model. Thus the comparison between our asymptotic variance and that obtained in the i.i.d. context and with the classical Hill estimator still remains valid. Concerning the estimator proposed by \cite{dHMZ2016}, they have obtained the asymptotic variance
$$\sigma^2(\theta, \gamma, \rho):= {\gamma^2 \over \rho^2} \left\{\left[(1-\rho)^2+\rho^2\right]r(1, 1)+2\rho(1-\rho){\theta^{1 \over \gamma} \log \theta^{1 \over \gamma} \over 1+\theta^{1\over \gamma}}\right\},$$
which is clearly again larger than our asymptotic variance ${\cal AV}(K_{\Delta^*_{opt}})$.

 \subsection{\textsf{Generalized autoregressive conditional heteroskedasticity (GARCH) model}}
  
We consider the GARCH model defined as
\begin{equation}
X_t=\sigma_t \, \varepsilon_t
\label{garch11}
\end{equation}
where $(\sigma_t)$ is a function of the history up to time $t-1$ represented by ${\cal{H}}_{t-1}$. The process of innovations $\varepsilon_t$ is a strict white noise with mean zero and
variance one and is assumed to be independent of ${\cal{H}}_{t-1}$. In other words $\sigma_t$ is ${\cal{H}}_{t-1}$-measurable,  ${\cal{H}}_{t-1}$ being the filtration generated by $(X_s)_{s\leq t-1}$ and therefore var$(X_t\mid {\cal{H}}_{t-1})=\sigma_t^2$.  The sequence $(X_t)$ follows a GARCH($p,q$) process if, for all $t$,
\begin{equation}
\sigma_{t}^2=\alpha_0+\sum_{j=1}^{p}\alpha_{j} X_{t-j}^{2}
+\sum_{k=1}^{q}\beta_{k}\sigma_{t-k}^2,\quad \alpha_j,\beta_k>0 .
\label{garchsigma}
\end{equation}
This model also satisfies  the regularity conditions $(C_R)$ but with a covariance structure $r$ which cannot be explicitly computed. In that case the comparison between the different estimators can be done only by simulation.

In fact, in Section \ref{simul}, we compare, in addition to the GARCH model, all the estimators for the three abovementioned models, through a simulation study. Actually, to be completely honest in the comparison, it is not sufficient to compare the constant in the variance. Rather, it is necessary to take the intermediate sequence into account as explained at the end of Section \ref{sectiongamma}. Indeed, the variance of our estimator is ${\cal AV}(K_{\Delta^*_{opt}})/k$ whereas that of \cite{dHMZ2016} is
$\sigma^2(\theta, \gamma, \rho)/\widetilde k$ with $\widetilde k$ of a larger order than $k$ due to the different conditions imposed.

\subsection{\textsf{Simulation study}}
\label{simul}

\begin{figure}[h!]
\begin{center}
\includegraphics[width=4cm]{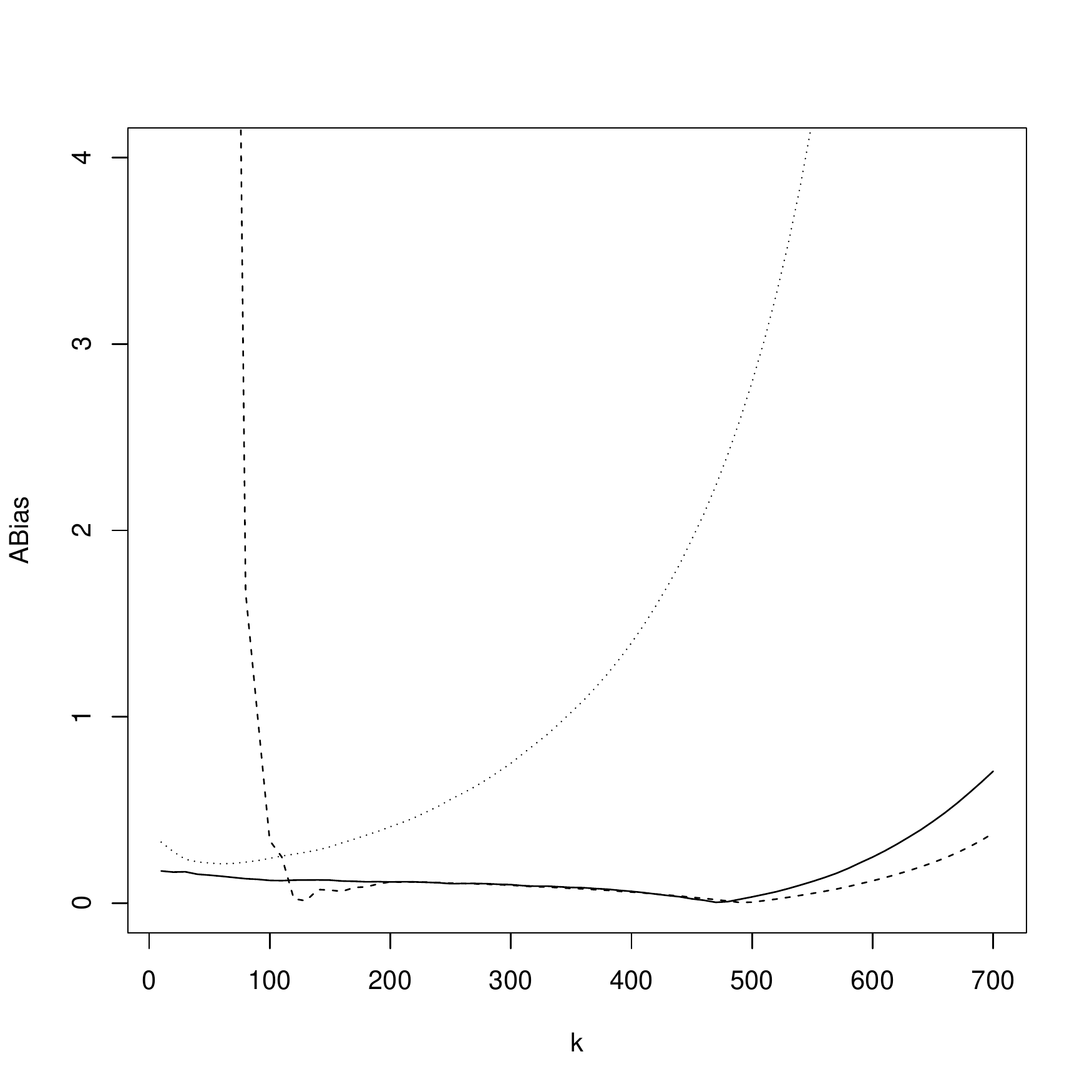}\includegraphics[width=4cm]{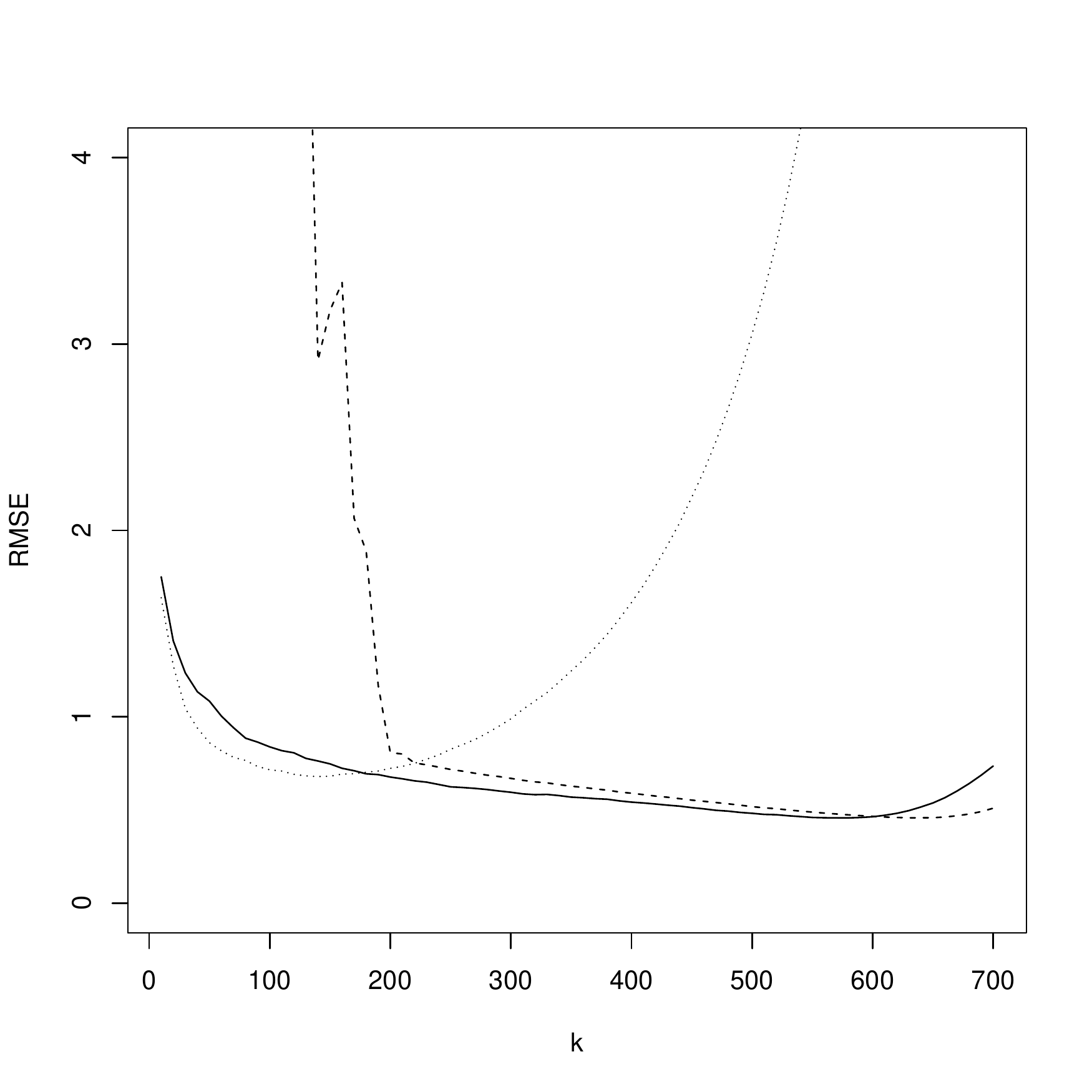}\\
\includegraphics[width=4cm]{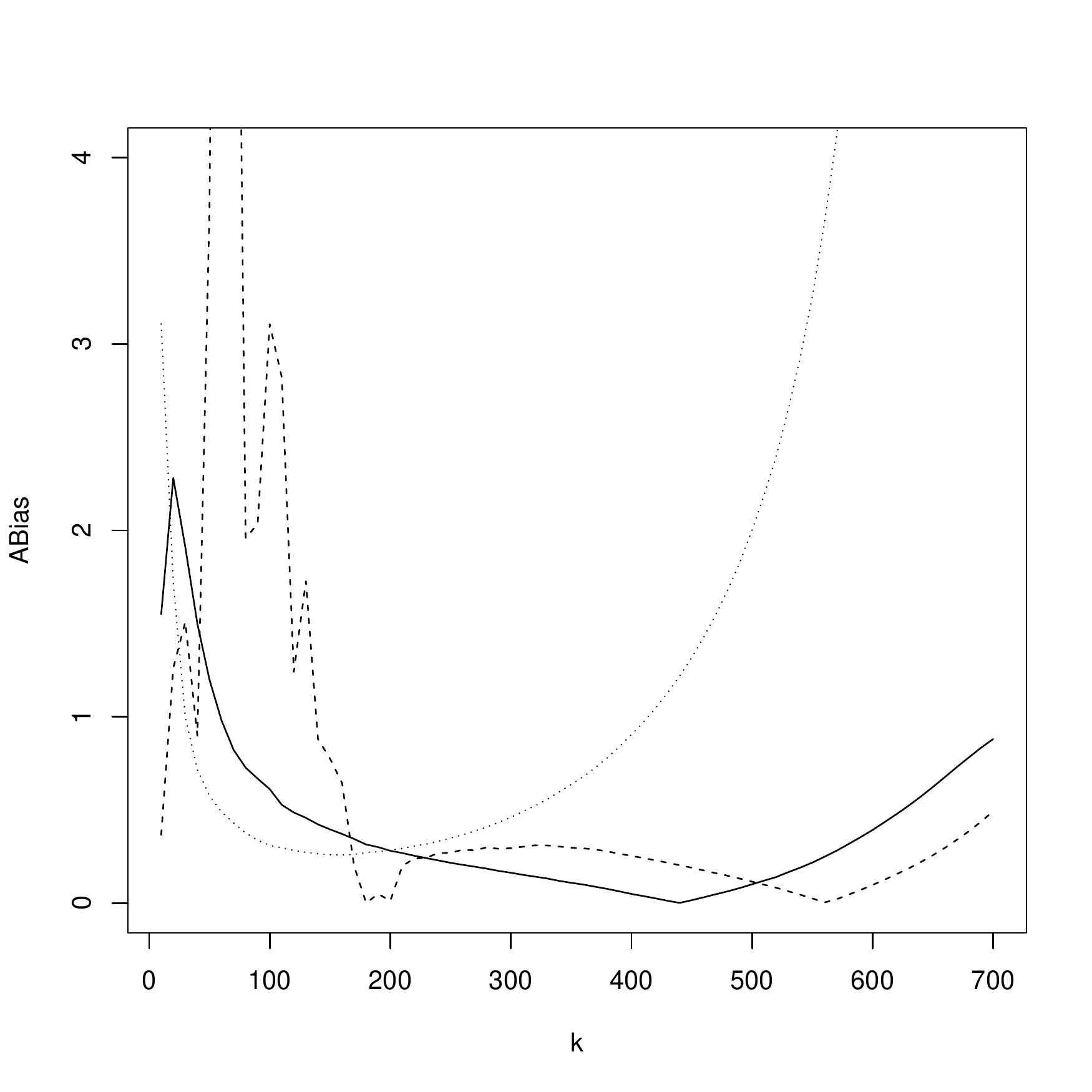}\includegraphics[width=4cm]{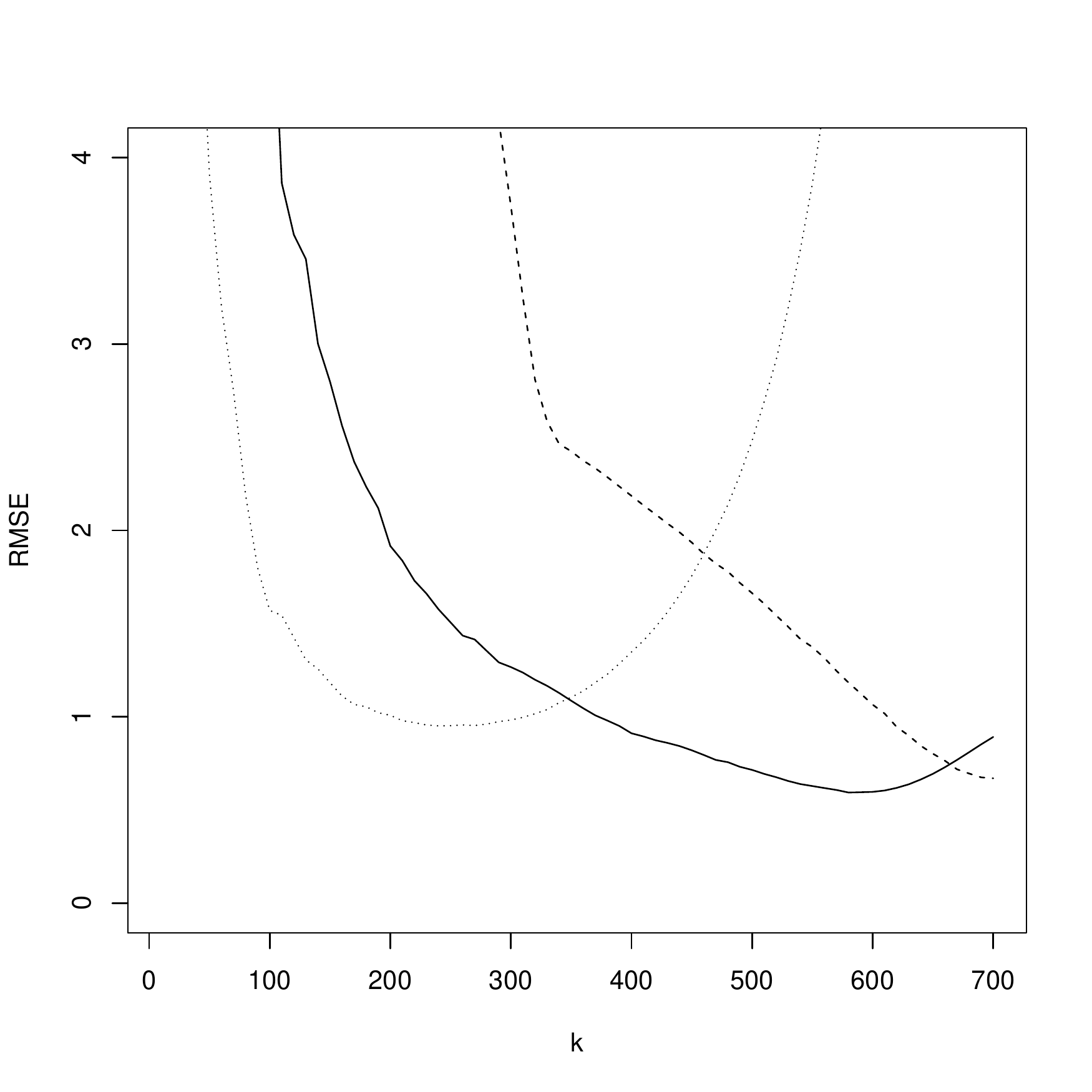}\\
\includegraphics[width=4cm]{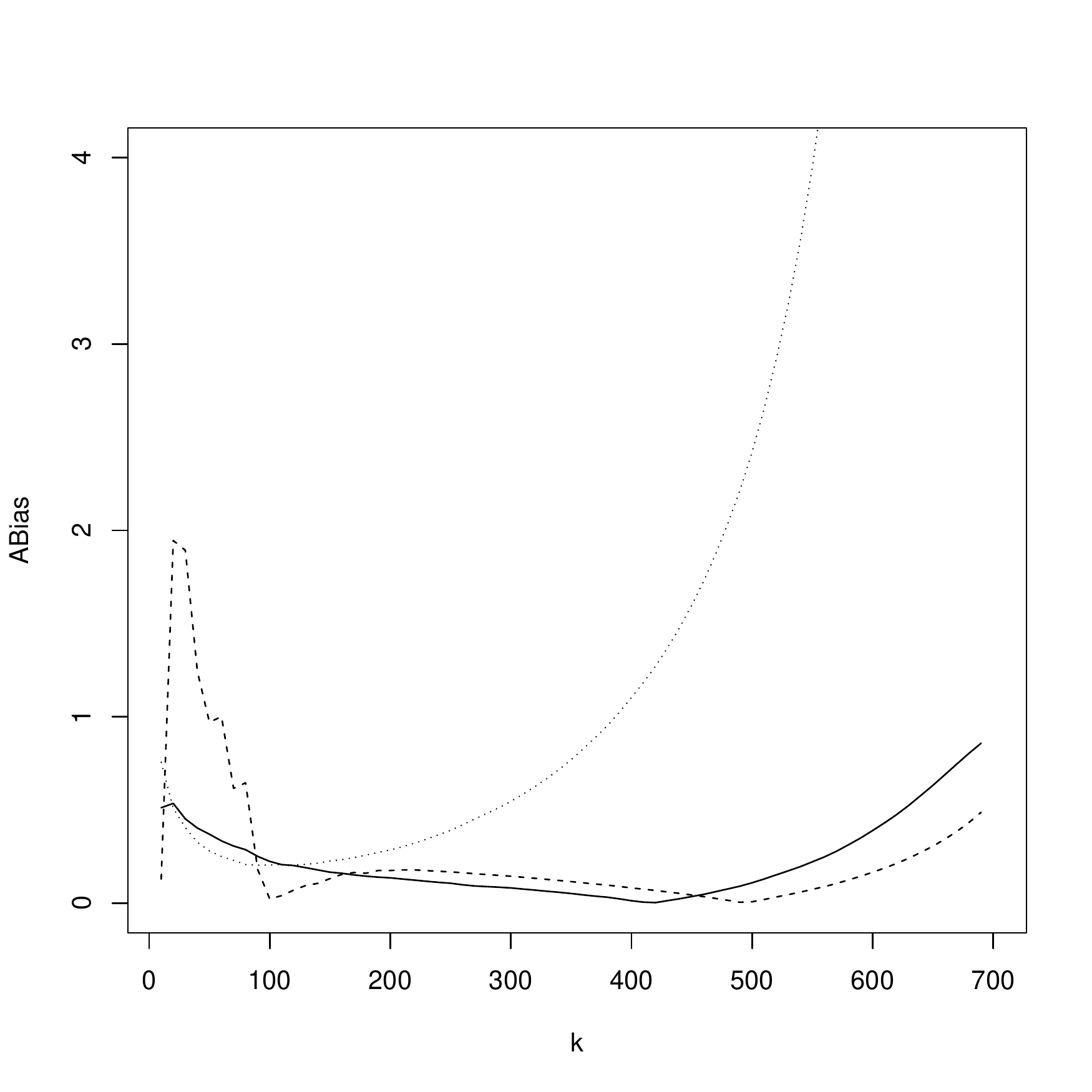}\includegraphics[width=4cm]{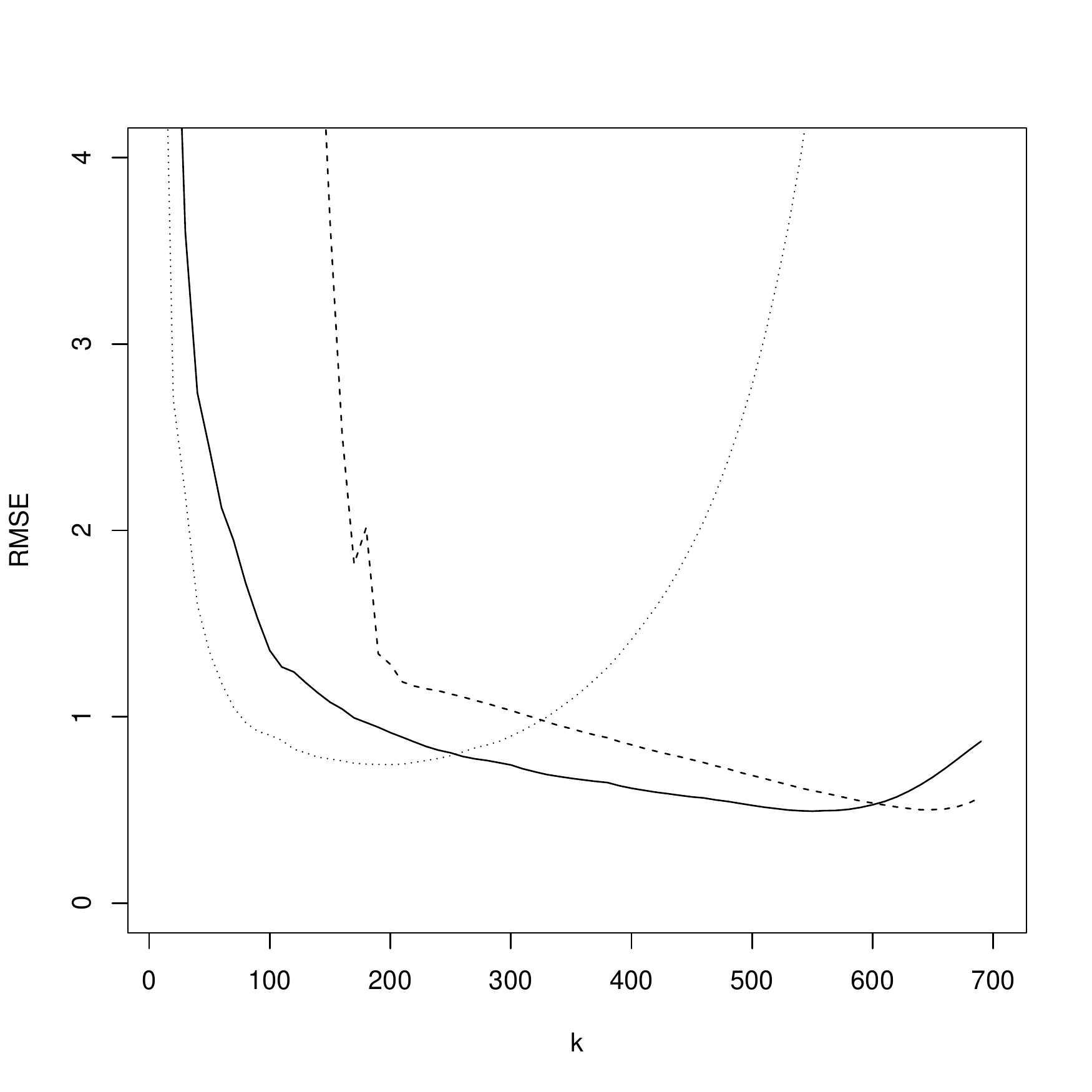}\\
\includegraphics[width=4cm]{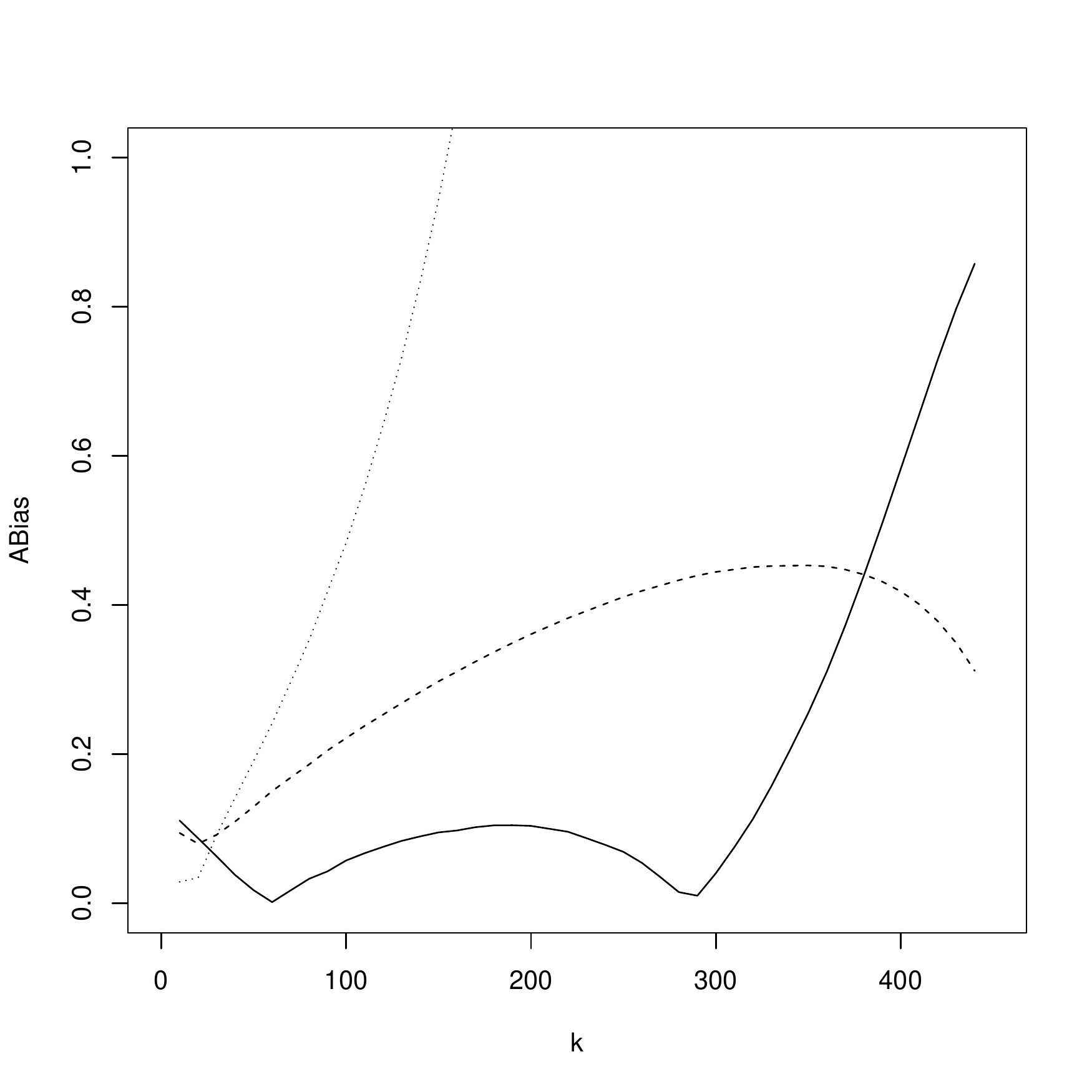}\includegraphics[width=4cm]{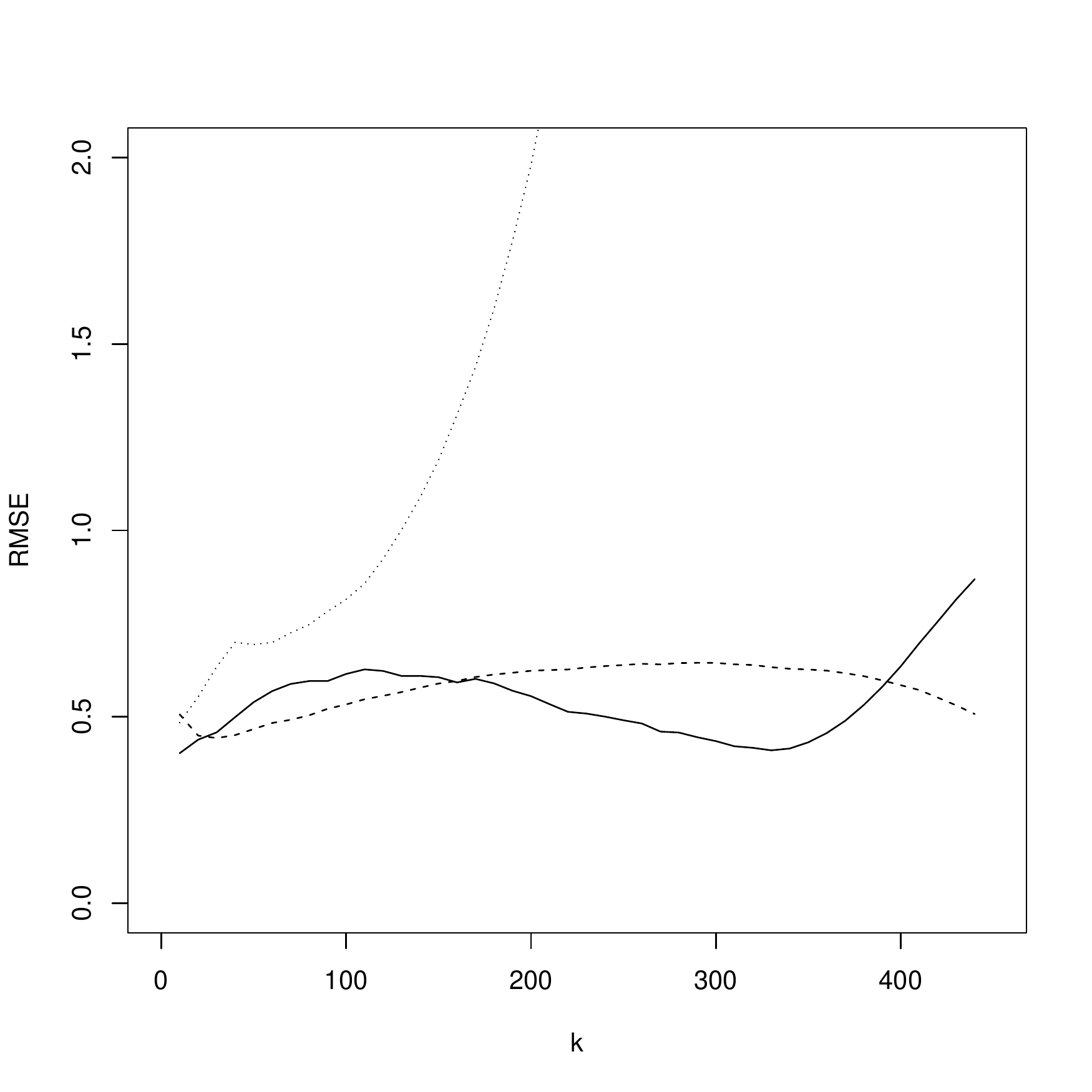}\\
\includegraphics[width=4cm]{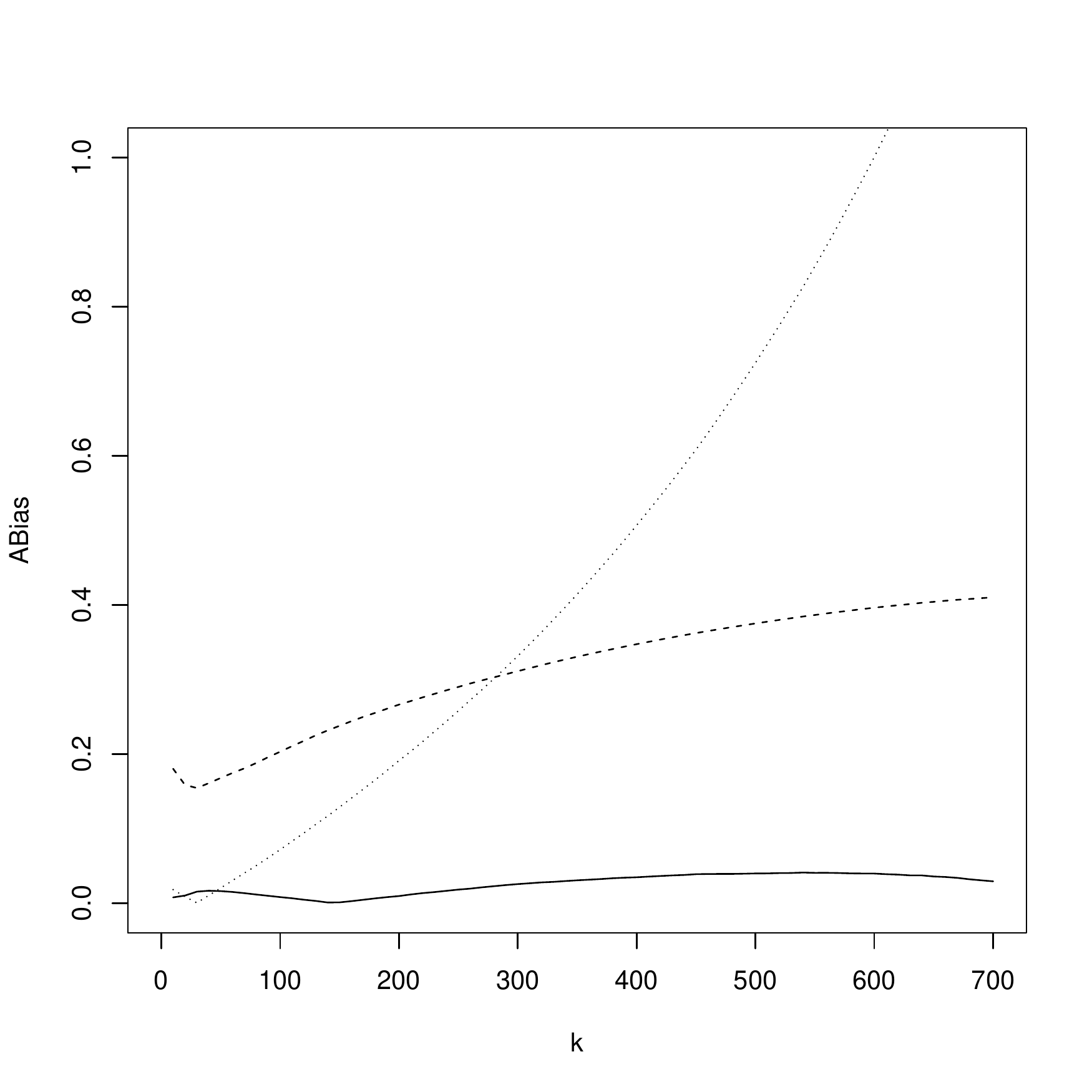}\includegraphics[width=4cm]{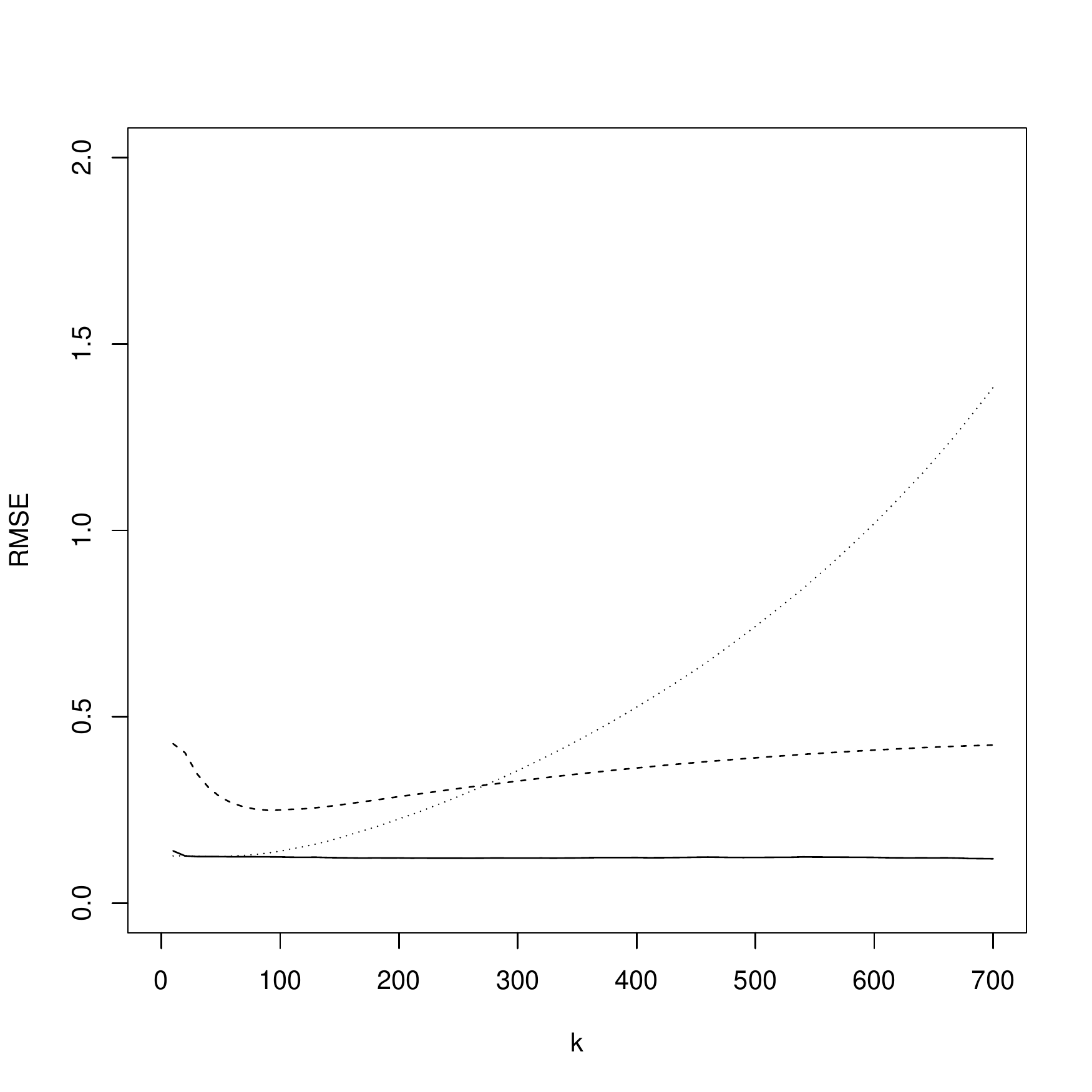}
\caption{\footnotesize Simulation study: By row, Models 1, 2, 3, 4, 5. By column, ABias (left) and RMSE (right). The full line shows our estimator, the dashed line the \cite{dHMZ2016} estimator and the dotted line represents the \cite{weiss1978} estimator.}
\label{Modelquant}
\end{center}
\end{figure}
We proceed to a simulation study to assess our high quantile estimator (\ref{quantileestimator}) with 
$p=0.001$  in five different model cases for which we can simulate the theoretical value of the true 99.9\% quantile. The three first models are the independence, AR(1) and MA(1) models proposed by \cite{dHMZ2016} and with $\varepsilon$ following the distribution
\begin{equation*}
F_\varepsilon(\varepsilon)=\left\{\begin{array}{ccc}(1-q)(1-\widetilde{F}(-\varepsilon)) & \rm{if} & \varepsilon<0,\\
1-q+q\widetilde{F}(\varepsilon)& \rm{if} & \varepsilon>0,\end{array}
\right.
\end{equation*}
where $\widetilde{F}$ is the unit Fr\'echet distribution function, $q=0.75$ so that $F_\varepsilon$ belongs to the max-domain of attraction with an extreme value index $\gamma=1$. We generate $N=5000$ time series of size $n=1000$ based on i.i.d  observations generated from $F_\varepsilon$ and we construct series from the three models
\begin{itemize}
\item Model 1: Independence model $X_i=\varepsilon_i$. The theoretical value of $x_{0.001}$ is 749.80.
\item Model 2: AR(1) model (\ref{AR}) with $\theta=0.3$. The theoretical value of $x_{0.001}$ is 1072.26.
\item Model 3: MA(1) model (\ref{MA}) with $\theta=0.3$. The theoretical value of $x_{0.001}$ is 972.85.
\end{itemize}
Note that the theoretical values are computed by Monte Carlo based on 1000 samples of size $10^6$.\\

We also consider two further models that are GARCH models with realistic parameters provided from our two real cases studied in Section \ref{applications}. They are 
\begin{itemize}
\item Model 4: GARCH(1,1) model (\ref{garch11}) with standardized Student $t$ innovations with 5.99 degrees of freedom and $\widehat{\alpha}_0=4.49\cdot 10^{-06}$, $\widehat{\alpha}_1=0.195$, $\widehat{\beta}_1=0.746$. The theoretical value of $x_{0.001}$ is 0.049. 
\item Model 5: GARCH(1,2) model (\ref{garch11}) with standardized Student $t$ innovations with 5.66 degrees of freedom and $\widehat{\alpha}_0=0.0443$, $\widehat{\alpha}_1=0.202$, $\widehat{\beta}_1=0.213$ and $\widehat{\beta}_2=0.467$. The theoretical value of $x_{0.001}$ is 3.103.
\end{itemize}
In both GARCH cases, the innovations being Student $t$, the distribution of $X_t$ belongs to the max-domain of attraction with an extreme value index evaluated at 0.27 for Model 4 and 0.15 for Model 5.
We simulate $N=5000$ time series of size $n=1000$ for Model 4 (corresponding to the sample size of our real financial data) and $n=4000$ for Model 5 (corresponding to the size of the wind speed data considered in our application). \\

Our high quantile (\ref{quantileestimator}) is based on our extreme value index estimator $\widehat \gamma_k(K_{\widehat\Delta^*_{opt}})$ and on the consistent estimator $\widehat \rho_{k_\rho}$ defined in (\ref{rhoestimateur}). Concerning the sequence $k_\rho$, we select it as follows
$$k_\rho:=\sup\left\{k: k\leq \min\left(m-1, {2m \over \log\log m}\right) \mbox{  and  } \widehat \rho_{k} \mbox{  exists}\right\},$$
with $m$ being the number of positive observations in the sample.
We compare our extreme quantile estimator $\widehat x_{p,\widehat\rho_{k_\rho}}$ with the one proposed by \cite{dHMZ2016} and defined as
\begin{eqnarray}
\widehat x_{k,k_\rho}(p)=X_{n-k,n} \left({k \over np}\right)^{\widehat \gamma_{k,k_\rho}} \left(1-{\left[M^{(2)}_k-2(\widehat \gamma^H_k)^2\right]\left[1-\widehat \rho_{k_\rho}\right]^2 \over 2 \widehat \gamma^H_k \widehat \rho^2_{k_\rho}}\left[1-\left({k\over np}\right)^{\widehat \rho_{k_\rho}}\right]\right),
\label{estimcorrige}
\end{eqnarray}
where $\widehat \rho_{k_\rho}$ is again the consistent estimator (\ref{rhoestimateur}) and $\widehat \gamma_{k,k_\rho}$ is their new estimator of the index defined as
$$
\widehat \gamma_{k,k_\rho}=\widehat \gamma^H_k-{\left[M^{(2)}_k-2(\widehat \gamma^H_k)^2\right]\left[1-\widehat \rho_{k_\rho}\right] \over 2 \widehat \gamma^H_k \widehat \rho_{k_\rho}}.
$$
Note that the estimator $\widehat x_{k,k_\rho}(p)$ is in fact slightly different from that included in the latter paper since a personal discussion with the authors allowed to identify an error in the proof of Theorem 4.2 in \cite{dHMZ2016}, precisely in their Assertion (A.6).\\

To compare the different estimators, we compute the absolute value of the mean of the bias (ABias) together with the root mean squared errors (RMSE) based on the $N$ samples, and defined as
\begin{eqnarray*}
\mbox{ABias}(x_p; k):=\left|{1 \over N}\sum_{i=1}^N {\widehat{x}_{p}^{(i)} \over x_p} - 1\right| \mbox{  and  } \mbox{RMSE}(x_p; k):= \sqrt{{1 \over N}\sum_{i=1}^N \left({\widehat{x}_{p}^{(i)} \over x_p} - 1\right)^2},
\end{eqnarray*}
where $\widehat{x}_{p}^{(i)}$ is the $i$-th value ($i=1,\ldots,N$) of the estimator of $x_p$ evaluated at $k$.
Figure \ref{Modelquant} shows the results for each of the five models by row and by column, the ABias (left) and RMSE (right). The full line corresponds to our high quantile estimator, the dotted line to the original \cite{weiss1978} estimator and the dashed line to the \cite{dHMZ2016} estimator. Note that it is not appropriate to compare our curves with that of \cite{dHMZ2016}, Figures 6,7,8, because our curves in Figure \ref{Modelquant} are computed from their corrected estimator. 
Globally, our high quantile estimator shows a lower bias than those of \cite{dHMZ2016} and Weissman especially for the GARCH models. The bias is also less variable than the two others for the lowest values of $k$. In terms of RMSE it is also very competitive when compared to the two alternative estimators. Both the bias and RMSE of our high quantile show a period of stability for a wider range of $k$; an important feature for practical applications.

\section{\textsf{REAL DATA ANALYSIS}}
\label{applications}
In this section, we illustrate the use of our estimator to calculate the daily VaR of a financial index series and the hourly return level of  wind speed data. Both the VaR and return level are a high quantile (\ref{quantile}) defined for a certain $p$-level. We use out-of-sample backtesting to assess the efficiency of our high quantile estimation. 
\subsection{\textsf{Financial index data}}
\begin{figure}[H]
\begin{center}
\includegraphics[width=12cm]{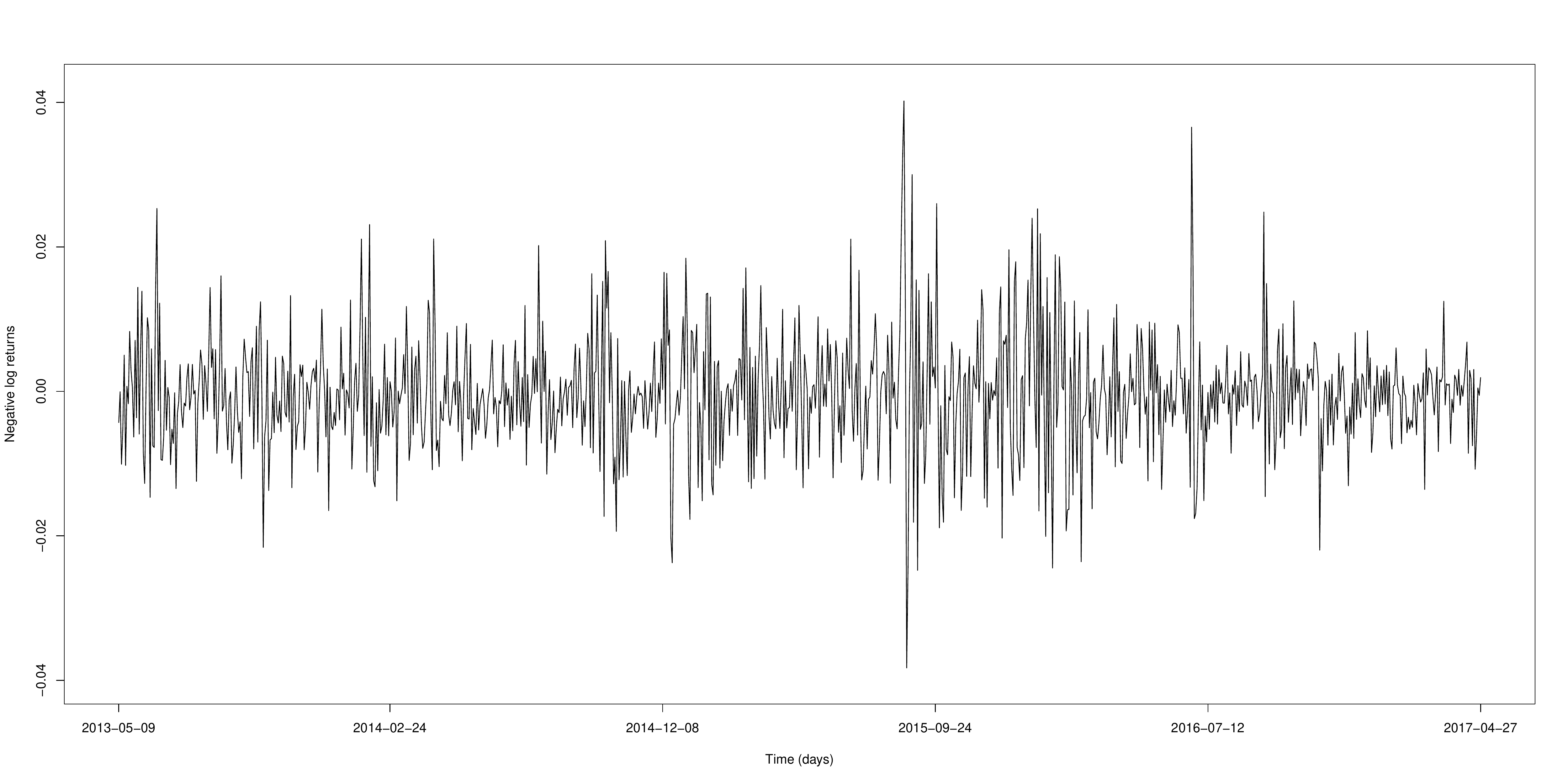}
\caption{\footnotesize S\&P500 index data: daily negative log-returns from 2013-05-09 to 2017-04-27.}
\label{SP}
\end{center}
\end{figure}

\begin{figure}[H]
\begin{center}
\includegraphics[width=6cm]{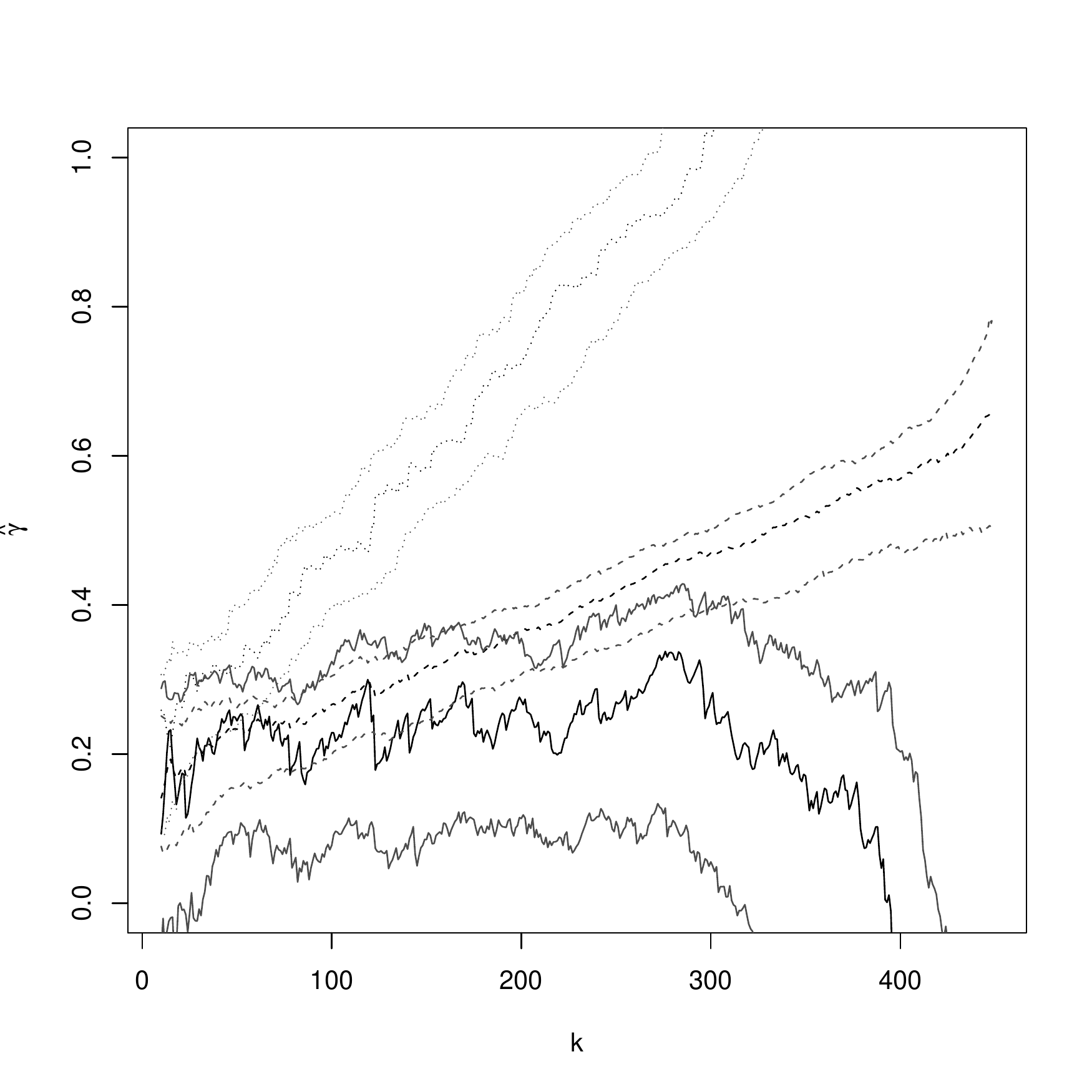}\includegraphics[width=6cm]{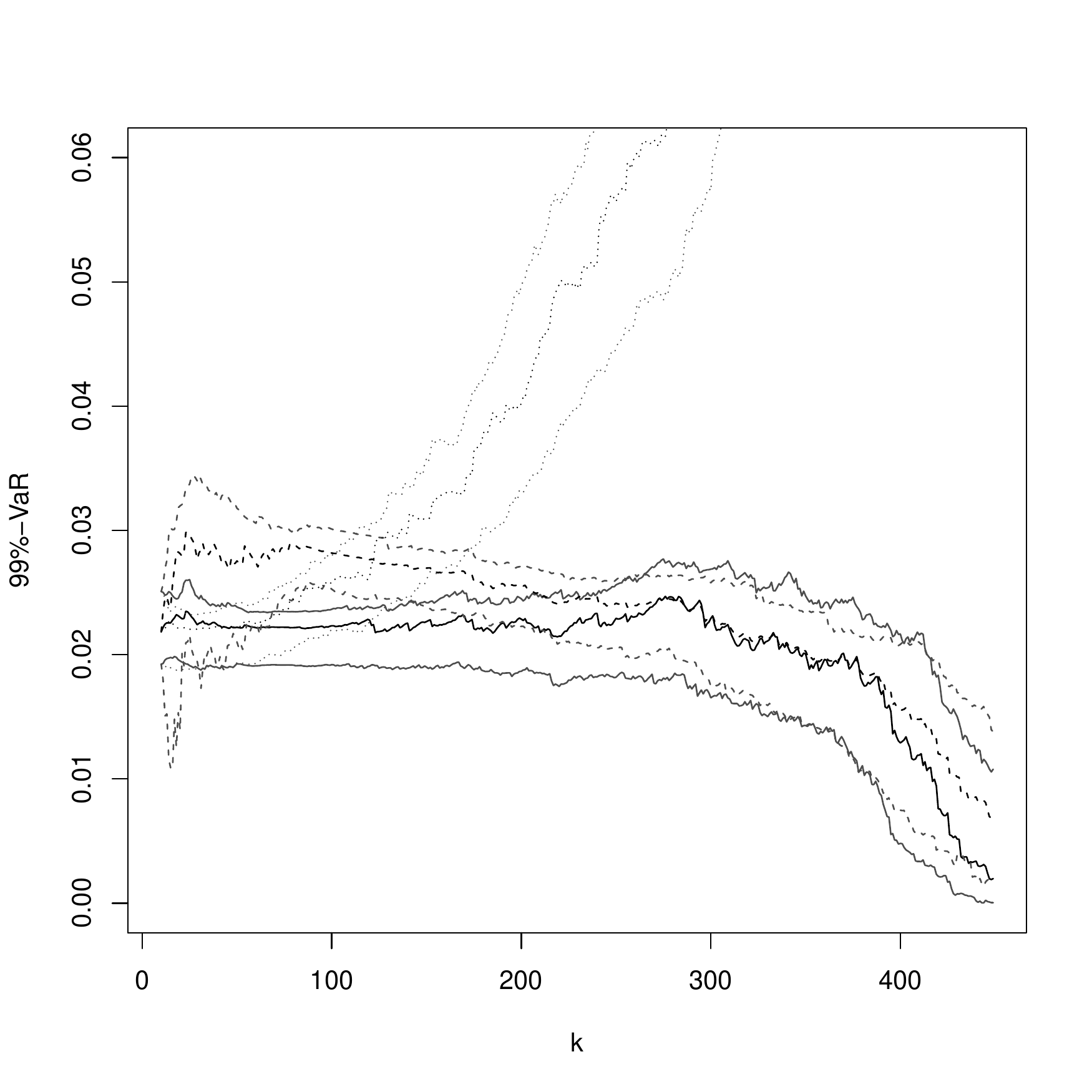}\includegraphics[width=6cm]{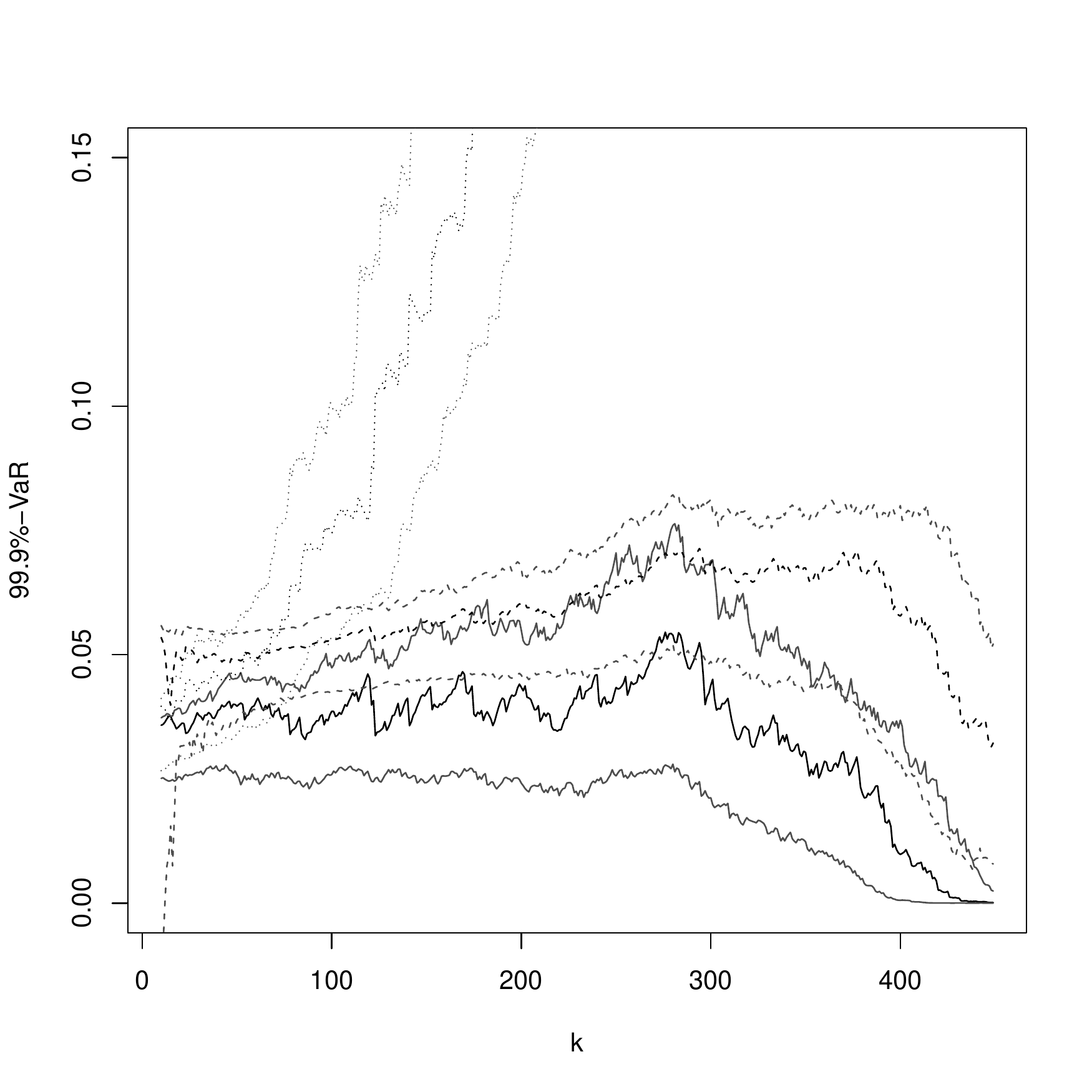}
\caption{\footnotesize S\&P500 index data: estimated values of $\gamma$ (left panel), 99\%-VaR (middle panel) and 99.9\%-VaR (right panel) against $k$ using our estimator $\widehat \gamma_k(K_{\widehat\Delta^*_{opt}})$ (left panel, full line), the \cite{dHMZ2016} estimator (dashed line), the Hill estimator (left panel, dotted line) and the Weissman estimator (middle and right panels, dotted line). The grey lines are the 95\%-bootstrap confidence intervals.}
\label{gammaSP}
\end{center}
\end{figure}

The data in Figure \ref{SP} shows the daily negative log-returns $X_t$ for $n=1000$ values of S\&P500 index from 2013-05-09 to 2017-04-27. In a risk management perspective, the Value-at-Risk (VaR) is a common quantity hardwired in
the international regulatory framework referred to as the Basel Accords \citep[see][for a historical review of the Basel International Settlement]{tarullo2008}. The Basel Accord requires the largest international banks to hold regulatory capital for the trading book based on a 99\%-VaR over a 1-day or 10-day holding period. The VaR-based risk capital calculation has received much attention over the last two decades \citep[see][Chapter 1]{mcneil2015}.  The $\alpha$\%-VaR for the horizon $h = 1$ day is the quantile $x_p$ (\ref{quantile}) with $p=1-\alpha$ of the distribution for the index daily log-returns. 
Stylized features of financial  series such as the S\&P500 index returns are heavy taildness; gaussianity assumption strongly violated; presence of heteroskedasticity or volatility clustering; absence of autocorrelations in returns \citep[see][for more details]{cont2005}. To model phenomenon with such characteristics,  the GARCH(1,1) model (\ref{garch11}) and (\ref{garchsigma}) with $p=q=1$ is specifically appropriate; see \cite{mcneil2015}, Chapter 4 for a review of GARCH models and Example 5.59, p. 235 for the GARCH fitting of financial time series. We fit a GARCH(1,1) model to our dataset with Student-$t$ innovations. The estimated model is 
 $$\sigma^2_t=\widehat{\alpha}_0+\widehat{\alpha}_1X^2_{t-1}+\widehat{\beta}_1\sigma^2_{t-1},$$
 with $\widehat{\alpha}_0=4.49\cdot 10^{-06}(1.34\cdot 10^{-06})$, $\widehat{\alpha}_1=0.195(0.0396)$, $\widehat{\beta}_1=0.746(0.0440)$ and the parameter of the Student $t$ is $\widehat{\nu}=5.99 (1.143)$ where the value in parentheses is the standard deviation. 
 The S\&P500 log-returns being of a stationary $\beta$-mixing type, we can use our estimator to calculate the $\alpha$\%-VaR for the horizon $h = 1$ day. We start to estimate the extreme value index $\gamma$ of the loss returns using our estimator $\widehat \gamma_k(K_{\widehat\Delta^*_{opt}})$ that we compare to the Hill estimator and the asymptotically unbiased estimator of \cite{dHMZ2016} for $\gamma$. Figure \ref{gammaSP}, left panel, shows the different curves of the estimated values against $k$. 
Our estimator (full line), even if more variable, looks more stable than the two others. Both the Hill estimator (dotted line) and that of \cite{dHMZ2016} (dashed line) increase with $k$. The gray lines are 95\%-confidence intervals calculated using block bootstrap \citep{buhlmann2002, davisonHinkley1997}. We use blocks of length 200 as suggested by \cite{dHMZ2016} and we simulate 99 bootstrap samples. We therefore estimate the 99\% (resp. 99.9\%)-VaR shown in Figure \ref{gammaSP} middle panel (resp. right panel) using our high quantile estimator (full line), \cite{dHMZ2016} estimator (dashed line) and \cite{weiss1978} estimator (dotted line). Whereas the Weissman estimator is not stable for both levels of $\alpha$, our estimator seems more stable than the one of \cite{dHMZ2016} which slightly decreases over $k$ for the lowest value of $\alpha$ and increases over $k$ for the highest $\alpha$.
The stability of the estimator is essential to decide a value of $k$ which will be used to get the high quantile estimator. The selection of $k$ is equivalent to the choice of the threshold in the EVT peaks-over-threshold method.  We arbitrarily choose $k=80$ that is a value within the window of stable values of $k$ in Figure \ref{gammaSP}. 
\begin{figure}[H]
\begin{center}
\includegraphics[width=8cm]{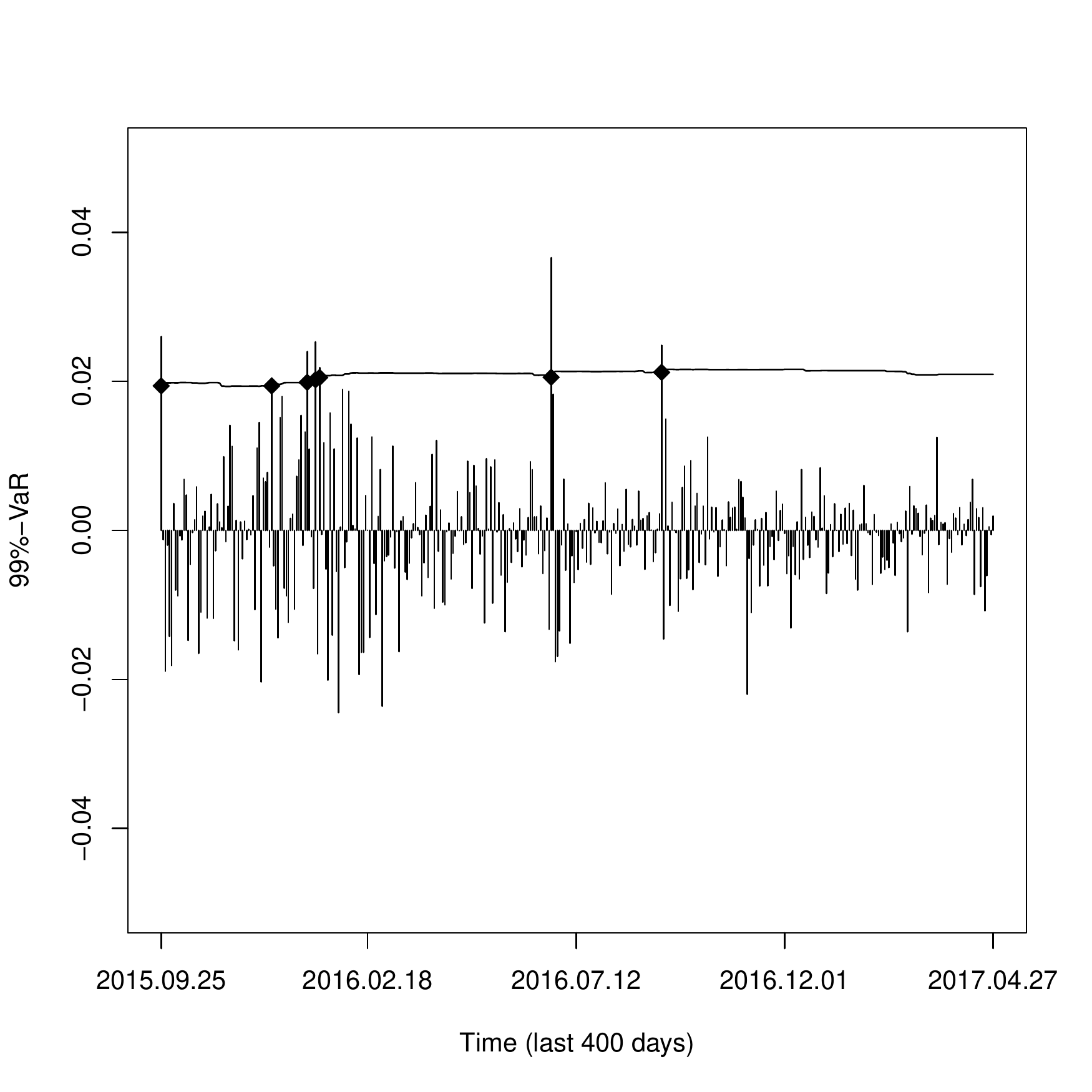}
\caption{\footnotesize S\&P500 index data: out-of-sample 99\%-VaR estimation for the last 400 daily returns of the observed period. The full line shows our high quantile estimate which uses at each time point a window of 600 historical data. The points are the times at which the realized log-return is higher than the estimated 99\%-VaR.}
\label{VaR99BacktestSP}
\end{center}
\end{figure}
To assess the efficiency of the VaR estimator, it is standard to use the so-called backtesting procedure \citep{danielsson2011} as suggested by the Basel Committee on banking supervision. The backtesting procedure consists of calculating the
$\alpha$-\%VaR estimate and comparing it to the value realized the next day. A violation
is said to occur whenever the estimated VaR is lower than the realized negative log-return. We repeat this operation over a window of historical data corresponding to the period 2015-09-24 to 2017-04-26 forecasting that way the last 400 daily 99\%-VaR. Figure \ref{VaR99BacktestSP} shows the resulting out-of-sample 99\%-VaR estimation calculated at each time point using a window of 600 past data. The expected number of violations of the VaR is 4 and the observed number is 7. The VaR violations are represented by the points in Figure \ref{VaR99BacktestSP}. The highest violation during this period happened on 2016-06-23 explained by the Wall Street reaction to Brexit (see, for instance, \url{http://www.reuters.com/article/us-usa-stocks-idUSKCN0Z918E}). We use the \cite{kupiec1995} coverage test that is a variation of the Binomial test and obtain a $p-$value of 0.173, failing to reject the null hypothesis of correct violation number. 

\subsection{\textsf{Wind speed data}}
\begin{figure}[h]
\begin{center}
\includegraphics[width=12cm]{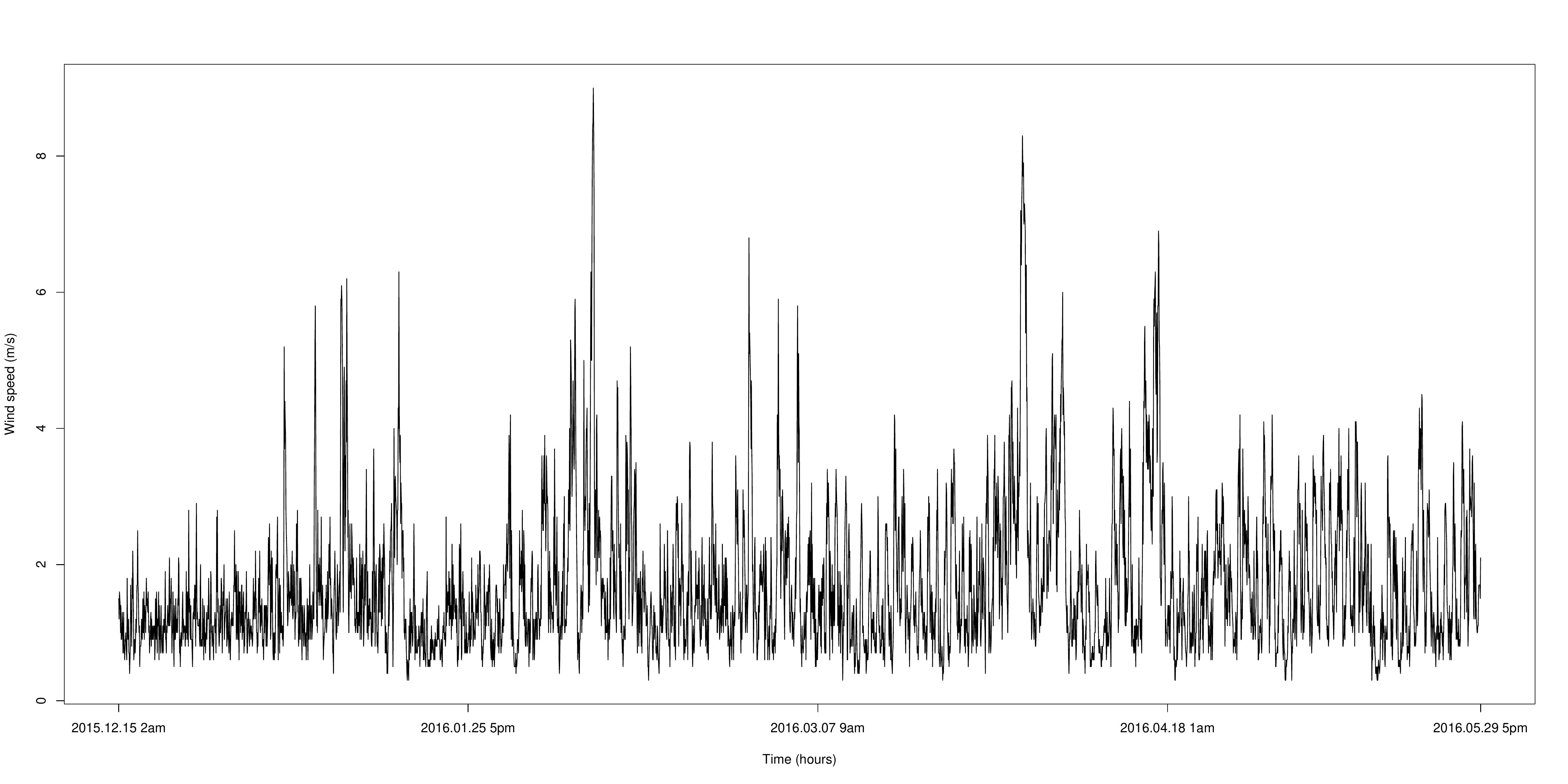}
\caption{\footnotesize Wind speed data: hourly measurements in Arosa, Switzerland, from 2015.12.15, 2am to 2016.05.29, 5pm.}
\label{vent}
\end{center}
\end{figure}
\begin{figure}[h]
\begin{center}
\includegraphics[width=8cm]{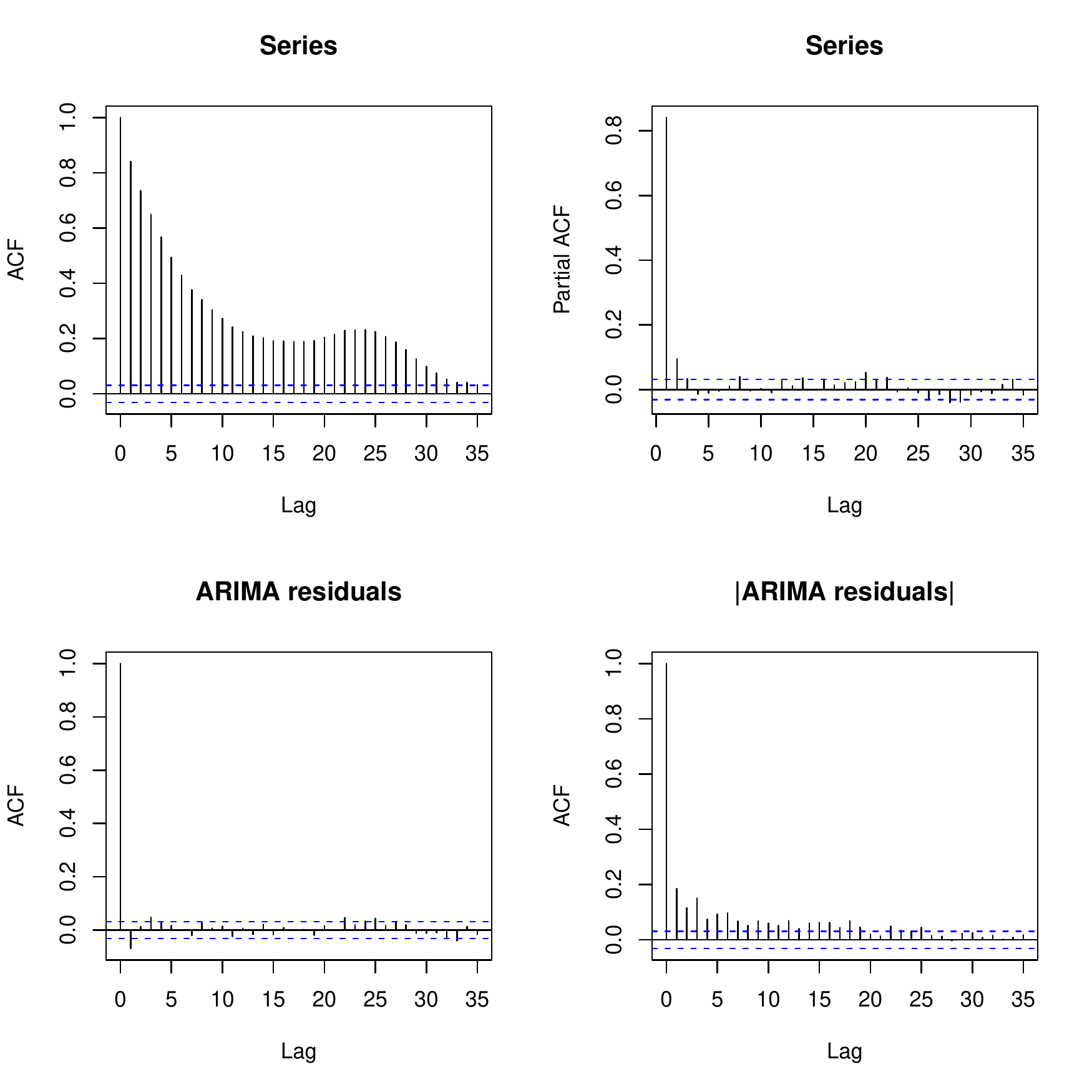}
\caption{\footnotesize Wind speed data: ACF (top left) and PACF (top right) of the wind speed series and ACF of the ARIMA residuals (bottom left) and of the absolute value of the ARIMA residuals (bottom right).}
\label{acf}
\end{center}
\end{figure}
\begin{figure}[h]
\begin{center}
\includegraphics[width=6cm]{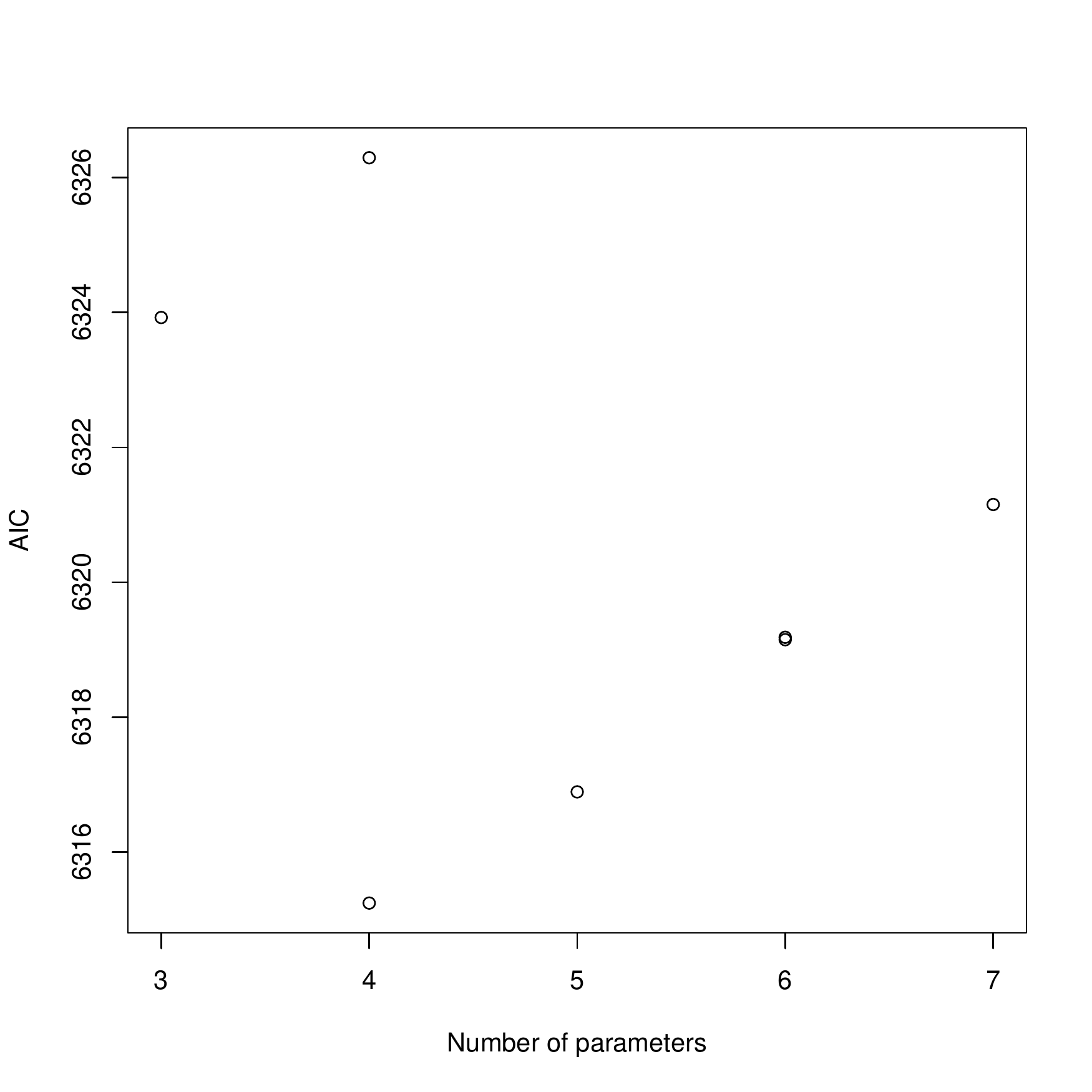}
\caption{\footnotesize Wind speed data: AIC against the GARCH number of parameters estimated.}
\label{aic}
\end{center}
\end{figure}
\begin{figure}[h]
\begin{center}
\includegraphics[width=6cm]{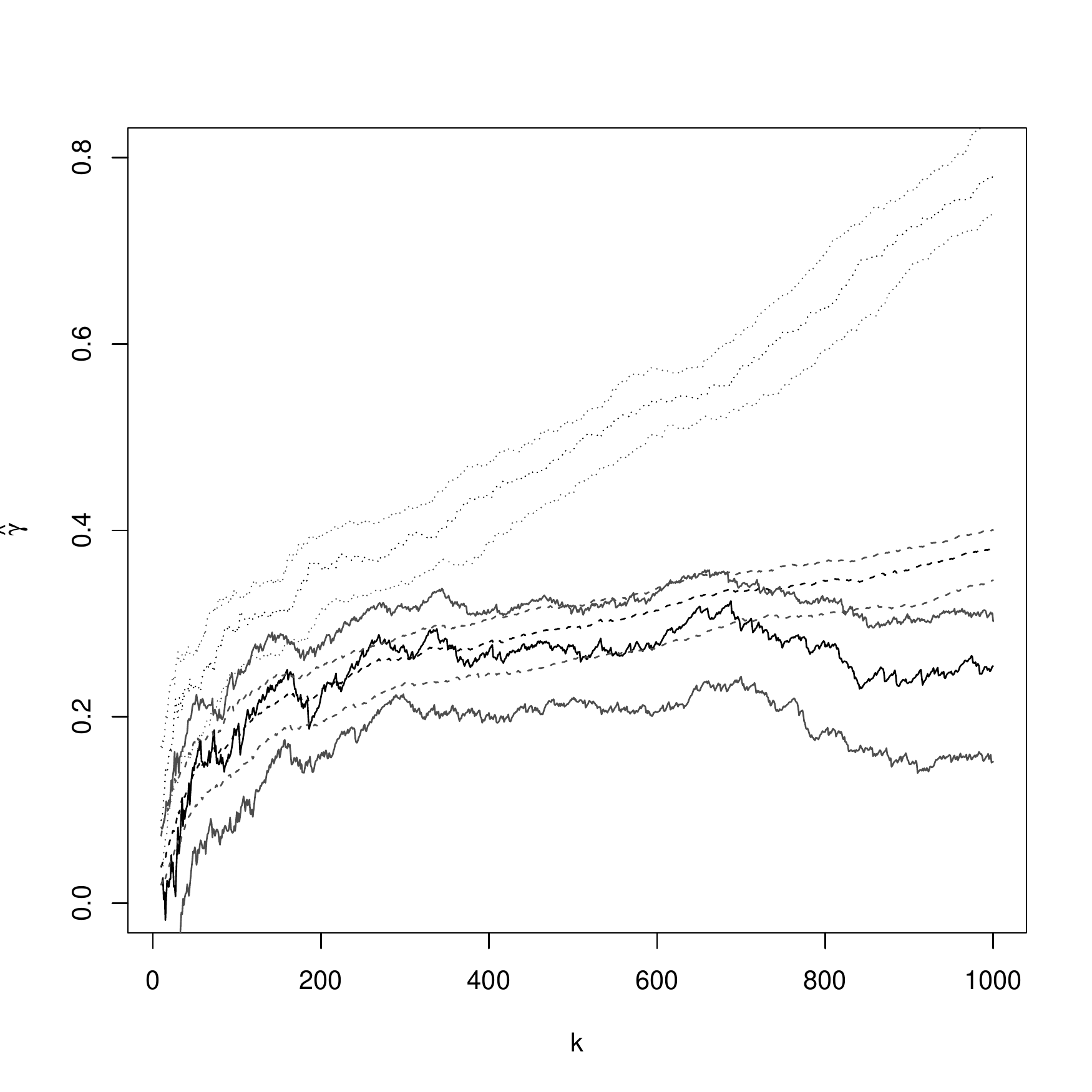}\includegraphics[width=6cm]{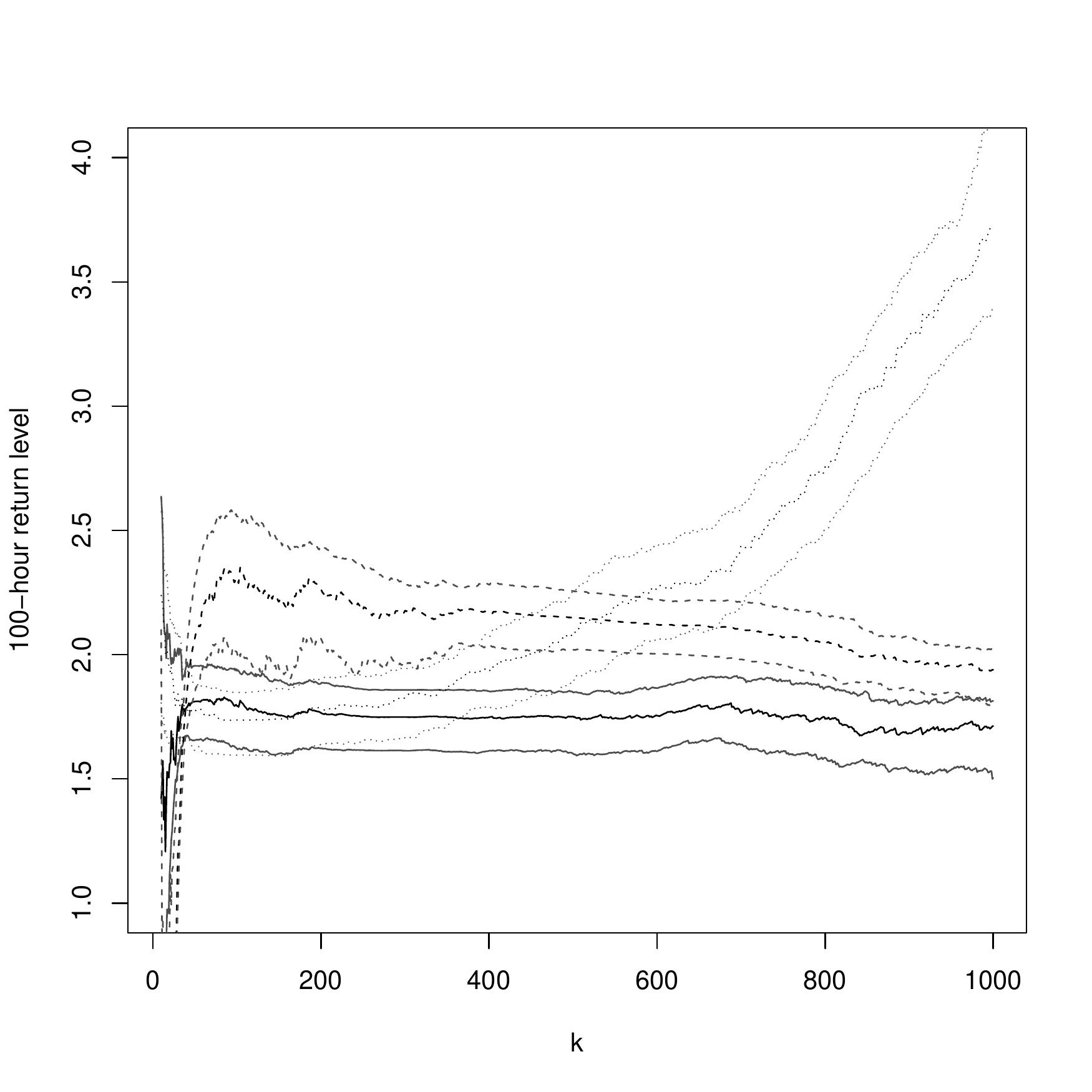}\includegraphics[width=6cm]{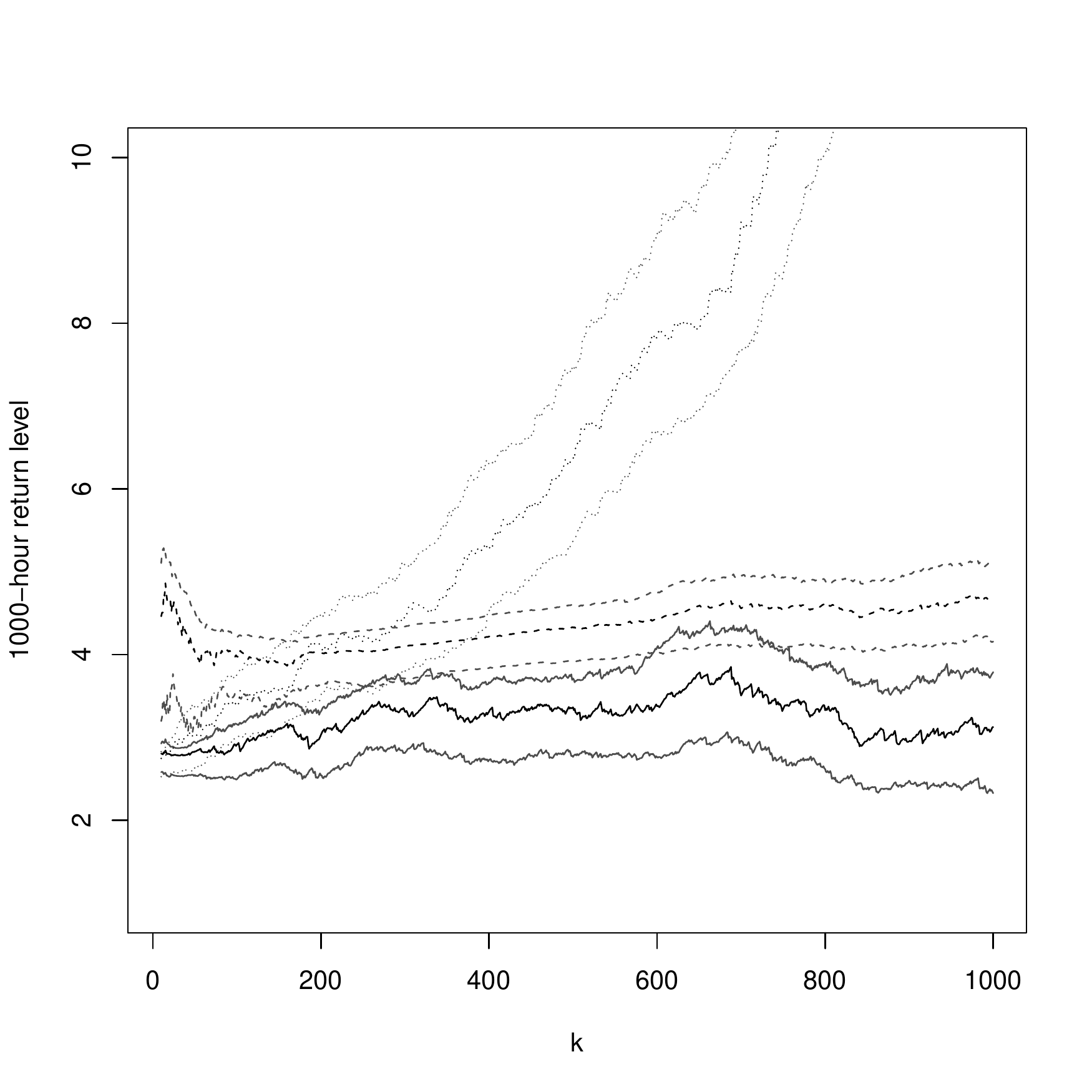}
\caption{\footnotesize Wind speed data: estimated values of $\gamma$ (left panel), 100-hour return level (middle panel) and 1000-hour return level (right panel) against $k$ using our estimator $\widehat \gamma_k(K_{\widehat\Delta^*_{opt}})$ (left panel, full line), the \cite{dHMZ2016} estimator (dashed line), the Hill estimator (left panel, dotted line) and the Weissman estimator (middle and right panels, dotted line). The grey lines are the 95\%-bootstrap confidence intervals.}
\label{gamma}
\end{center}
\end{figure}

Figure \ref{vent} shows hourly wind speed (m/s) data $X_t$, ($t$ represents hour) measured in Arosa, Switzerland, from 2015.12.15, 2am to 2016.05.29, 5pm, consisting of $n=3895$ hourly values. The data were provided by the Federal Office of Meteorology and Climatology, SwissMeteo. The data show evidence of seasonality confirmed by the ACF plot of the series in the top left panel of Figure \ref{acf}. Several recent papers suggest that wind speed data are of an ARMA-GARCH-type. \cite{Lojowska2010}, for instance, claim that artificial wind speeds simulated from ARMA-GARCH models are statistically indistinguishable from the real wind speed time series measurements under observation; \cite{LSE2013} use the ARMA-GARCH model to predict the time series mean and volatility of wind speed. The top left panel of Figure \ref{acf} shows that autocorrelations of Arosa wind speeds are significant for a large number of lags but we can notice that the PACF plot (top right panel) has a significant spike only at lag 1, meaning that all the higher-order autocorrelations are effectively explained by the lag-1 autocorrelation. More formally, we fit different ARIMA models and the model selected according to AIC is an ARIMA(1,1,1) corresponding to a model with one AR term and one MA term and a first difference used to account for a linear trend in the data as follows
\begin{eqnarray*}
Y_t=X_t-X_{t-1}\\
Y_t=\phi_1 Y_{t-1}+\varepsilon_t-\theta_1 \varepsilon_{t-1},
\end{eqnarray*}
where $\varepsilon_t$ is a random shock occurring at time $t$.  Fitted to the data, the estimated parameters are $\widehat{\phi}_1 =0.819(0.0125)$, $\widehat{\theta}_1=-0.989(0.0050)$.
The ACF in the bottom left panel of Figure \ref{acf} shows no serial correlation remaining for the ARIMA residuals
\begin{eqnarray}
e_t & = & x_t-x_{t-1}-\widehat{\phi}_1x_{t-1}+\widehat{\phi}_1x_{t-2}+\widehat{\theta}_1e_{t-1}, \hspace{1cm}t=1,\ldots,n.
\label{residuals}
\end{eqnarray}
However,  from the ACF plot of the absolute values for the ARIMA residuals $e_1,\ldots,e_n$, there is evidence of volatility clustering. The absence of autocorrelations for the residuals and the presence of volatility clustering clearly suggests a GARCH-type model. This is confirmed by a Ljung-Box test on the squared ARIMA residuals (\ref{residuals}) providing a $p-$value lower than $10^{-5}$, rejecting the null hypothesis that the series is a strict white noise.
We fit several GARCH models of different orders $(p,q)$ on the ARIMA residuals $e_1,\ldots,e_n$. The ARIMA residuals being heavy-tailed, we use Student $t$ innovations. The GARCH form is given by (\ref{garch11}) where here $X_t$ represents our residuals and $\sigma_t$ is the ARIMA wind speed residual volatility. We select the orders $p$ and $q$ using AIC. Figure \ref{aic} represents the AIC against the total number of parameters estimated for each model, that is $p+q+1$. The AIC is minimized for $p+q+1=4$. Not readable from the graph, the lowest AIC value for $p+q+1=4$ corresponds to the model GARCH(1,2) which is
 $$\sigma^2_t=\widehat{\alpha}_0+\widehat{\alpha}_1e^2_{t-1}+\widehat{\beta}_1\sigma^2_{t-1}+
 \widehat{\beta}_2\sigma^2_{t-2},$$
 $\widehat{\alpha}_0=0.0443(0.00946)$, $\widehat{\alpha}_1=0.202(0.0279)$, $\widehat{\beta}_1=0.213(0.0764)$ and $\widehat{\beta}_2=0.467(0.0787)$ with standardized Student $t$ innovations with degrees of freedom $\widehat{\nu}=5.66$. 
With such data, a measure of extreme events of interest is the $1/p$-hour return level $r(p)$ with $p$ small. The return level is the value that has a $p\%$ chance of being exceeded in a given hour. Our proposed estimator for high quantile (\ref{quantile}) can be used on the series of the ARIMA residuals $e_1,\ldots,e_n$ being GARCH-type removed from seasonality and therefore satisfying the stationary $\beta$-mixing conditions. Note that to obtain a return level for the wind speed original data $r_{x_t}(p)$ (in m/s) at time $t$ one can use the ARIMA model and the $1/p$-hour return level $r_{e_t}(p)$ estimate of the ARIMA residual $e_t$ at time $t$, that is
\begin{equation}
\widehat{r}_{x_t}(p)  = \widehat{r}_{e_t}(p) - \widehat{\theta}_1e_{t-1}+x_{t-1}+\widehat{\phi}_1x_{t-1}-\widehat{\phi}_1x_{t-2}.
\label{transfo}
\end{equation}
To estimate the high quantile of the ARIMA residuals  $r_{e}(p)$, we first estimate $\gamma$ using our estimator $\widehat \gamma_k(K_{\widehat\Delta^*_{opt}})$. The full line of the left panel in Figure \ref{gamma} shows our estimated values of $\gamma$ against different values of $k$. Compared to the Hill (dotted line) and the \cite{dHMZ2016} (dashed line) estimators, our estimator looks more stable over $k$ even if it seems more variable. The grey lines are the 95\%-confidence intervals calculated using a block bootstrapping method with block length of size 200 and based on 99 bootstrap samples. For the high quantile estimation with $p=0.01$ (resp. $p=0.001$), corresponding to the 100-hour return level (resp. 1000-hour return level), our estimator (full line) stays very stable as shown in the middle panel of Figure \ref{gamma} (resp. right panel) compared to the Weissman estimator (dotted line) and even to that of \cite{dHMZ2016} (dashed line) which shows a slightly decreasing trend (resp. increasing trend). 

\begin{figure}[h]
\begin{center}
\includegraphics[width=14cm]{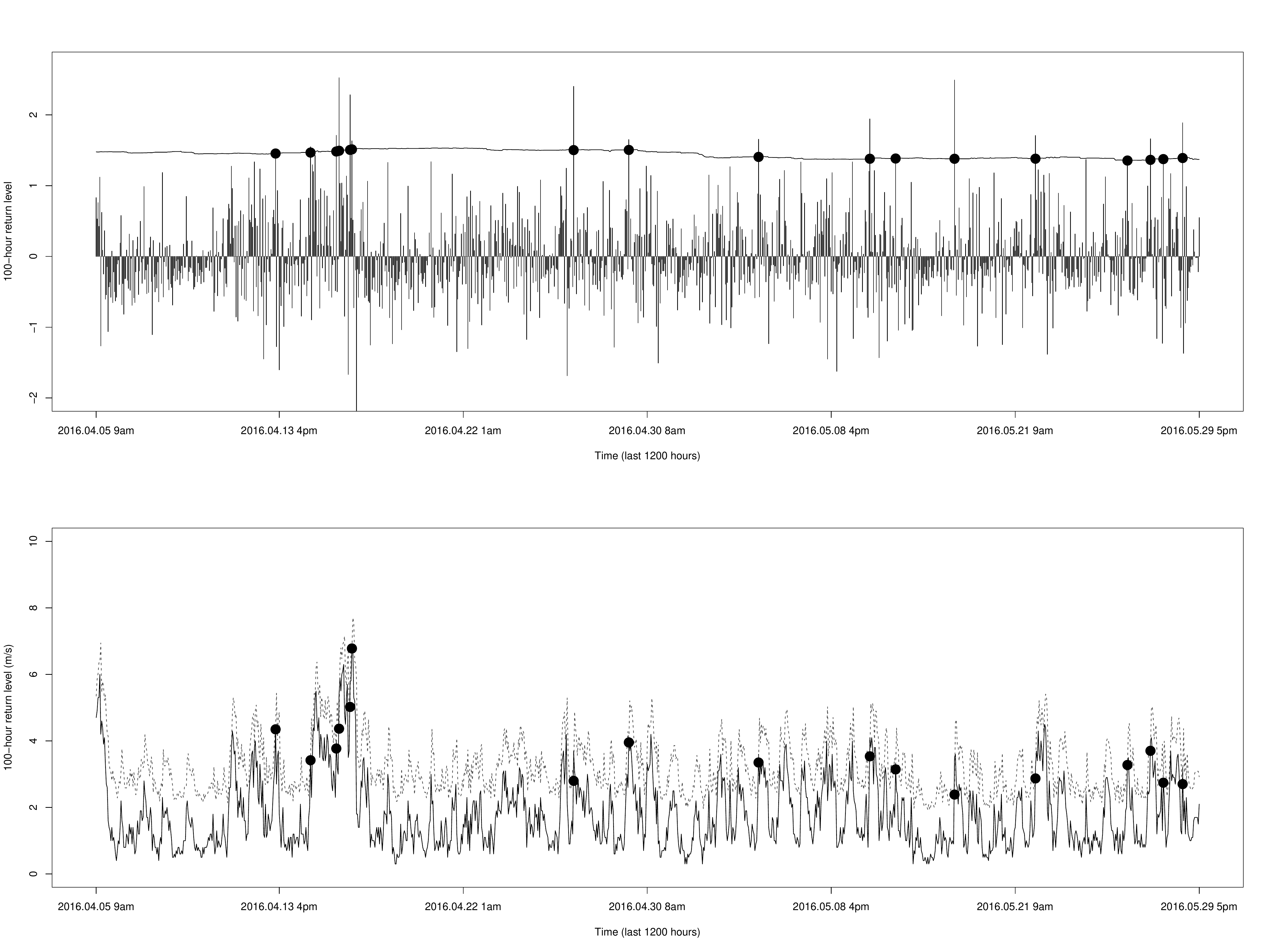}
\caption{\footnotesize Wind speed data: out-of-sample 100-hour return level estimation forecasting the last 1200 hourly returns of the observed period. The full line shows our high quantile estimate which uses at each time point a window of 2000 historical data. The top panel shows the high quantile estimate of the ARIMA residuals. The bottom panel shows the high quantile estimate of the wind speed data in m/s. The points are the times at which the ARIMA residual (top) and wind speed value (bottom) is higher than the estimated return level.}
\label{VaR999backtestVent}
\end{center}
\end{figure}

From Figure \ref{gamma}, we can reasonably choose any value of $k$ between 200 and 600 without leading much to variable results. We chose $k=450$ for our estimator and as for the financial data, we proceed to an out-of-sample estimation of the 100-hour return level. The line in Figure \ref{VaR999backtestVent} (top panel) shows our 100-hour return level estimates for the ARIMA residuals evaluated at each time point of the 1200 hours from 2016.04.05 at 8am to 2016.05.29 at 4pm forecasting the next hour 99\%-quantile till the last hour of the observed period. To estimate the return level at each point, we use a window of the 2000 past values. The expected number of 100-hour return level violations is 12 and the observed number is 17. The Kupiec coverage test for the correct number of exceedances is not rejected with a $p$-value of 0.172. The bottom panel of Figure~\ref{VaR999backtestVent} shows the 100-hour return level for the original wind speed data (in m/s) calculated using (\ref{transfo}). Confidence intervals are not shown in the plot for clarity sake. Because of the two-step method used (filtering the seasonality using ARIMA model and then applying our estimator on the residuals), a convenient way to get confidence intervals is by proceeding to a block bootstrap of the original data and applying the two-step method for each resample.

\section{\textsf{CONCLUSION}}
\label{conclusion}
In this paper we have introduced a new asymptotically unbiased estimator of high quantiles for $\beta$-mixing stationary time series. Comparing the new procedure to the alternative proposed by \cite{dHMZ2016}, our high quantile estimator provides, in addition to lower ABias and RMSE in general, more stability over $k$, an important feature expected in this type of approach to be applicable in practice. In application, the new high quantile estimator can be proposed to any other stationary $\beta$-mixing heavy-tailed time series for which high quantiles needed to be calculated. This concerns heavy-tailed autoregressive data encountered in network traffic forecasting for instance and many other applications data in climate change.

\clearpage
\section*{\textsf{APPENDIX: PROOFS OF THE RESULTS}}
\label{appendix}
Before establishing our Theorem~1, we need  a result similar to Proposition A.1 in \cite{dHMZ2016} but under the weak assumptions of our Theorem~1, excluding $(C_K)$. In particular, a third order condition is not assumed.

{\bf Proposition 1.} {\it Let $(X_1, X_2, ...)$ be a stationary $\beta$-mixing time series with a continuous common marginal distribution function $F$ and assume $(C_{S0})$ and $(C_R)$. Suppose that $k$ is an intermediate sequence such that $\sqrt k A(n/k)=O(1)$. For a given $\varepsilon >0$, under a Skorohod construction, there exist a function $\widetilde A\sim A$, and a centered Gaussian process $(W(t))_{t \in [0, 1]}$ with covariance function $r$, such that, as $n\to \infty$
$$\sup_{t\in(0, 1]} t^{{1\over 2}+\varepsilon}\left|\sqrt k\left(\log {Q_n(t)\over U({n\over k})} + \gamma \log t\right)-\gamma t^{-1}W(t)-\sqrt k \, \widetilde A\left({n \over k}\right) {t^{-\rho}-1\over \rho} \right|\longrightarrow 0 \quad {a.s.}$$
}
{\bf Proof of Proposition 1.} It is similar to that of Proposition A.1 in \cite{dHMZ2016} but assuming that $\sqrt k A\left({n \over k}\right)=O(1)$ and without a third order condition, thus below we only give the main differences.  In our context, the key inequality is the following: for all $\varepsilon, \delta>0$, there exists some positive number $u_0=u_0(\varepsilon,\delta)$ such that for $ux\geq u_0$:
\begin{eqnarray}
\left|{\log U(ux)-\log U(u)-\gamma \log x \over \widetilde A(u)} - {x^\rho-1 \over \rho} \right| \leq \varepsilon x^\rho \max\left(x^\delta, x^{-\delta}\right),
\label{inequality}
\end{eqnarray}
see, for instance, Theorem B.2.18 in \cite{dHF2006}. Now, using the representation $X_i=U(Y_i)$ where $Y_i$ follows a standard Pareto distribution, $(Y_1, Y_2, ...)$ is a stationary $\beta$-mixing series satisfying the regularity conditions $(C_R)$. Then, according to \cite{drees2003}, since $Q_n(t)=U(Y_{n-[kt],n})$ and under a Skorohod construction, there exists a centered Gaussian process $(W(t))_{t\in [0, 1]}$ with a covariance function $r$ such that for $\varepsilon>0$, as $n\to \infty$
$$\sup_{t\in (0, 1]} t^{{1\over 2}+\varepsilon} \left|\sqrt k \left(t {Y_{n-[kt], n} \over n/k}-1\right)-t^{-1} W(t)\right|\longrightarrow 0 \quad {a.s.}$$
This convergence combining with inequality (\ref{inequality}) entails that
$$\left|\log Q_n(t)-\log U\left({n\over k}\right)-\gamma\log\left({k\over n} Y_{n-[kt],n}\right) - \widetilde A\left({n\over k}\right) {\left({k\over n} Y_{n-[kt],n}\right)^\rho-1\over \rho}\right|\leq \varepsilon \left|\widetilde A\left({n\over k}\right)\right| \left({k\over n} Y_{n-[kt],n}\right)^{\rho+\delta}$$
for sufficiently large $n>n_0(\varepsilon, \delta)$, with probability 1.

Consequently
\begin{eqnarray*}
t^{{1\over 2}+\varepsilon}\left|\sqrt k\left(\log {Q_n(t)\over U(n/k)} + \gamma \log t\right)-\gamma t^{-1}W(t)-\sqrt k \, \widetilde A\left({n \over k}\right) {t^{-\rho}-1\over \rho} + \sqrt k \widetilde A\left({n \over k}\right){1\over \rho} \left\{t^{-\rho} - \left({k\over n} Y_{n-[kt],n}\right)^\rho\right\}\right.\\
\left.-\gamma\left\{\sqrt k \left(\log\left({k\over n} Y_{n-[kt],n}\right)+\log t\right)-t^{-1}W(t)\right\}
\right|\leq \varepsilon \sqrt k \left|\widetilde A\left({n\over k}\right)\right| t^{{1\over 2}+\varepsilon}   \left({k\over n} Y_{n-[kt],n}\right)^{\rho+\delta}.
\end{eqnarray*}
Choosing $\delta\in (0, -\rho)$, Proposition 1 then follows similarly as Proposition A.1 in \cite{dHMZ2016} since $\varepsilon$ can be arbitrarily close to 0. $\blacksquare$

{\bf Proof of Theorem~1.} From Proposition 1, we can easily infer that
\begin{eqnarray*}
\sqrt k \left\{\int_0^1 \log {Q_n(t) \over Q_n(1)} d(tK(t))+\gamma \int_0^1 \log t \, d(tK(t))\right\} = \gamma\int_0^1 \left[t^{-1}W(t)-W(1)\right] d(tK(t))\\
+\sqrt k \widetilde A\left({n\over k}\right)\int_0^1 {t^{-\rho}-1\over \rho} \, d(tK(t))+o(1)\int_0^1 t^{-{1\over 2}-\varepsilon} \, d(tK(t)).
\end{eqnarray*}
Using integration by parts, we have
\begin{eqnarray*}
\int_0^1 \log t \, d(tK(t))&=&-1,\\
\mbox{and }\int_0^1 {t^{-\rho}-1\over \rho} \, d(tK(t))&=&\int_0^1 t^{-\rho}K(t)dt, 
\end{eqnarray*}
our Theorem~1 now follows under $(C_K)$ by taking $0<\varepsilon<{1\over 2}-\tau$. $\blacksquare$

{\bf Proof of Corollary 1.} It is a direct consequence of Theorem~1 since $(W(t))_{t\in [0, 1]}$ is a centered Gaussian process with covariance function $r$. $\blacksquare$

{\bf Proof of Corollary 2.} It is a direct consequence of Corollary 1 since by construction the bias of $\widehat \gamma_k(K_{\Delta^*})$ is null. $\blacksquare$

{\bf Proof of Corollary 3.} According to Corollary 1, we only need to check the bias term. Recall that
$$K_{\widetilde \Delta^*_{opt}}(t)=\left({1-\widetilde\rho\over \widetilde\rho}\right)^2 - {(1-\widetilde\rho)(1-2\widetilde\rho)\over \widetilde\rho^2} t^{-\widetilde\rho},$$
from which we deduce that
$$\int_0^1 t^{-\rho} K_{\widetilde \Delta^*_{opt}}(t) dt =  {(1-\widetilde \rho)(\widetilde \rho-\rho) \over \widetilde \rho(1-\rho)(1-\widetilde\rho-\rho)}.$$
This achieves the proof of Corollary 3. $\blacksquare$

{\bf Proof of Theorem~2.}  Let $\Delta^*_{opt}:=\left({1-\rho \over \rho}\right)^2$ and $\widehat \Delta^*_{opt}:=\left({1-\widehat\rho_{k_\rho} \over \widehat\rho_{k_\rho}}\right)^2$. We consider the decomposition
\begin{eqnarray}
\sqrt k \left(\widehat \gamma_k(K_{\widehat \Delta^*_{opt}})-\gamma\right)&=&
\sqrt k \left(\widehat \gamma_k(K_{\Delta^*_{opt}})-\gamma\right)
+ \sqrt k \left(\widehat \gamma_k(K_{\widehat \Delta^*_{opt}})-\widehat \gamma_k(K_{\Delta^*_{opt}})\right).
\label{decompositionnewbis}
\end{eqnarray}
According to Corollary 2 we have
$$\sqrt k \left(\widehat \gamma_k(K_{\Delta^*_{opt}})-\gamma\right)
\stackrel{d}{\longrightarrow} {\cal N}\left(0, {\cal{AV}}(K_{\Delta^*_{opt}})\right).$$
To prove Theorem~2, it is thus sufficient to show that the second term in (\ref{decompositionnewbis}) is $o_\mathbb P(1)$. For this aim, note that
\vskip2ex
\noindent
$\displaystyle{
\sqrt k \left(\widehat \gamma_k(K_{\widehat \Delta^*_{opt}})-\widehat \gamma_k(K_{\Delta^*_{opt}})\right)}$
\begin{eqnarray*}
&=&\sqrt k \left\{\int_0^1 \log {Q_n(t) \over Q_n(1)} d(tK_{\widehat \Delta^*_{opt}}(t))-\int_0^1 \log {Q_n(t) \over Q_n(1)} d(tK_{ \Delta^*_{opt}}(t))\right\}\\
&=&\sqrt k \left\{\widehat \Delta^*_{opt}\int_0^1 \log {Q_n(t) \over Q_n(1)} dt-\Delta^*_{opt} \int_0^1 \log {Q_n(t) \over Q_n(1)} dt\right.\\
&&\left. +(1-\widehat \Delta^*_{opt})\int_0^1 \log {Q_n(t) \over Q_n(1)} d(tK_{2,\widehat\rho_{k_\rho}}(t))-(1-\Delta^*_{opt}) \int_0^1 \log {Q_n(t) \over Q_n(1)} d(tK_{2,\rho}(t))\right\}\\
&=&\sqrt k \left(\widehat \Delta^*_{opt}-\Delta^*_{opt}\right)\left\{\int_0^1 \log {Q_n(t) \over Q_n(1)} dt-\int_0^1 \log {Q_n(t) \over Q_n(1)} d(tK_{2,\rho}(t))\right\}\\
&&+\sqrt k \left(1-\widehat \Delta^*_{opt}\right)\left\{\int_0^1 \log {Q_n(t) \over Q_n(1)} d(tK_{2,\widehat \rho_{k_\rho}}(t))-\int_0^1 \log {Q_n(t) \over Q_n(1)} d(tK_{2,\rho}(t))\right\}\\
&=&\left(\widehat \Delta^*_{opt}-\Delta^*_{opt}\right) \left\{\sqrt k \left(\widehat \gamma_k(K_1)-\gamma\right)-\sqrt k \left(\widehat \gamma_k(K_{2,\rho})-\gamma\right)\right\}\\
&&+\left(1-\widehat \Delta^*_{opt}\right)\, \sqrt k \, \left\{\int_0^1 \log {Q_n(t) \over Q_n(1)} d(tK_{2,\widehat \rho_{k_\rho}}(t))-\int_0^1 \log {Q_n(t) \over Q_n(1)} d(tK_{2,\rho}(t))\right\}\\
&=:& T_1+T_2.
\end{eqnarray*}
We will study the two terms separately.\\
{\bf Term $T_1$.}  Using the consistency in probability of $\widehat \rho_{k_\rho}$ and the convergences  
\begin{eqnarray*}
\sqrt k \left(\widehat \gamma_k(K_1)-\gamma\right)
&\stackrel{d}{\longrightarrow}& {\cal N}\left(\lambda/(1-\rho), {\cal{AV}}(K_1)\right)\\
\sqrt k \left(\widehat \gamma_k(K_{2,\rho})-\gamma\right)
&\stackrel{d}{\longrightarrow}& {\cal N}\left(\lambda {\cal {AB}}(K_{2,\rho}), {\cal{AV}}(K_{2,\rho})\right)
\end{eqnarray*}
coming from Corollary 1, we have $T_1=o_\mathbb P(1)$.
\vskip2ex
\noindent
{\bf Term $T_2$.} For $\varepsilon\in (0, 1/2)$, uniformly for $t\in (0, 1]$:
$$\log {Q_n(t) \over Q_n(1)}=\gamma(-\log t)+{\gamma\over \sqrt k} \left[t^{-1}W(t)-W(1)\right]+\widetilde A\left({n \over k}\right){t^{-\rho}-1\over \rho}+{o(1)\over \sqrt k} t^{-{1\over 2}-\varepsilon}.$$
This implies that
\vskip2ex
\noindent
$\displaystyle{ \sqrt k \left\{
\int_0^1 \log {Q_n(t) \over Q_n(1)} d(tK_{2,\widehat \rho_{k_\rho}}(t))-\int_0^1 \log {Q_n(t) \over Q_n(1)} d(tK_{2,\rho}(t))\right\}}$
\begin{eqnarray*}
&=&\gamma \, \sqrt k \left\{\int_0^1 (-\log t)d(tK_{2,\widehat\rho_{k_\rho}}(t)) -  \int_0^1 (-\log t)d(tK_{2,\rho}(t))\right\}\\
&&+\gamma \left\{\int_0^1  \left[t^{-1}W(t)-W(1)\right]d(tK_{2,\widehat\rho_{k_\rho}}(t)) - 
\int_0^1  \left[t^{-1}W(t)-W(1)\right]d(tK_{2,\rho}(t))\right\}\\
&&+\sqrt k \widetilde A\left({n \over k}\right) \left\{\int_0^1   {t^{-\rho}-1\over \rho}d(tK_{2,\widehat\rho_{k_\rho}}(t)) -
 \int_0^1   {t^{-\rho}-1\over \rho}d(tK_{2,\rho}(t))\right\}\\
&&+o(1) \int_0^1 t^{-{1\over 2}-\varepsilon} d(tK_{2,\widehat\rho_{k_\rho}}(t)) +o(1) \int_0^1 t^{-{1\over 2}-\varepsilon} d(tK_{2,\rho}(t))\\
&=&\gamma \, \sqrt k \left\{\left[-t\log t \left(K_{2,\widehat\rho_{k_\rho}}(t)-K_{2,\rho}(t)\right)\right]_0^1+\int_0^1 K_{2,\widehat\rho_{k_\rho}}(t)dt - \int_0^1 K_{2,\rho}(t)dt\right\}\\
&&+\gamma \left\{ \int_0^1\left[t^{-1}W(t)-W(1)\right]\left(K_{2,\widehat\rho_{k_\rho}}(t)-K_{2,\rho}(t)\right)dt
+\int_0^1\left[t^{-1}W(t)-W(1)\right] t \left(K'_{2,\widehat\rho_{k_\rho}}(t)-K'_{2,\rho}(t)\right)dt\right\}\\
&&+\sqrt k \widetilde A\left({n \over k}\right)\left\{\left[t\, {t^{-\rho}-1\over \rho}\left(K_{2,\widehat\rho_{k_\rho}}(t)-K_{2,\rho}(t)\right)\right]_0^1+\int_0^1 t^{-\rho} \left(K_{2,\widehat\rho_{k_\rho}}(t)-K_{2,\rho}(t)\right)dt\right\}\\
&&+o(1)\left\{\left[t^{{1\over 2}-\varepsilon}\left(K_{2,\widehat\rho_{k_\rho}}(t)-K_{2,\rho}(t)\right)\right]_0^1+\left({1\over 2}+\varepsilon\right)\int_0^1 t^{-{1\over 2}-\varepsilon} \left(K_{2,\widehat\rho_{k_\rho}}(t)-K_{2,\rho}(t)\right)dt\right\}\\
&&+o(1) \int_0^1 t^{-{1\over 2}-\varepsilon} d(tK_{2,\rho}(t))\\
&=&\gamma \left\{ \int_0^1\left[t^{-1}W(t)-W(1)\right]\left(K_{2,\widehat\rho_{k_\rho}}(t)-K_{2,\rho}(t)\right)dt
+\int_0^1\left[t^{-1}W(t)-W(1)\right]t \left(K'_{2,\widehat\rho_{k_\rho}}(t)-K'_{2,\rho}(t)\right)dt\right\}\\
&&+o_\mathbb P(1). 
\end{eqnarray*}
Now, note that for $\varepsilon\in (0, 1/4)$ and $\widetilde \rho$ a random value between $\rho$ and $\widehat\rho_{k_\rho}$, we have
\vskip2ex
\noindent
$\displaystyle{
\bullet \left| \int_0^1\left[t^{-1}W(t)-W(1)\right]\left[K_{2,\widehat\rho_{k_\rho}}(t)-K_{2,\rho}(t)\right]dt \right|}$
\begin{eqnarray*}
&\leq&  \int_0^1 \left|t^{-1}W(t)-W(1)\right|\left|K_{2,\widehat\rho_{k_\rho}}(t)-K_{2,\rho}(t) \right| dt\\
&\leq& (1-\widehat\rho_{k_\rho}) \int_0^1 \left|t^{-1}W(t)-W(1)\right| \left|t^{-\widehat\rho_{k_\rho}}-t^{-\rho}\right|dt+|\widehat\rho_{k_\rho}-\rho|\int_0^1 \left|t^{-1}W(t)-W(1)\right|dt\\
&\leq& (1-\widehat\rho_{k_\rho}) \sup_{t\in(0, 1]} \left|t^{{1\over 2}+\varepsilon}\left[t^{-1}W(t)-W(1)\right]\right| \sup_{t\in(0, 1]} t^{{1\over 4}} \left|t^{-\widehat\rho_{k_\rho}}-t^{-\rho}\right| \int_0^1 t^{-{3\over 4}-\varepsilon}dt\\
&&+|\widehat\rho_{k_\rho}-\rho|\int_0^1 \left|t^{-1}W(t)-W(1)\right|dt\\
&\leq& {4\over 1-4\varepsilon} \left| \widehat\rho_{k_\rho}-\rho \right|(1-\widehat\rho_{k_\rho}) \sup_{t\in(0, 1]} \left|t^{{1\over 2}+\varepsilon}\left[t^{-1}W(t)-W(1)\right]\right| \sup_{t\in(0, 1]} (-\log t)t^{{1\over 4}-\widetilde\rho} \\
&&+|\widehat\rho_{k_\rho}-\rho|\int_0^1 \left|t^{-1}W(t)-W(1)\right|dt\\
&=&o_\mathbb P(1),
\end{eqnarray*}
since $\sup_{t\in (0, 1]} t^{{1\over 2}+\varepsilon} t^{-1} |W(t)|=O(1)$ a.s.. Similarly, we have \\
\vskip2ex
\noindent
$\displaystyle{
\bullet \left| \int_0^1\left[t^{-1}W(t)-W(1)\right]t\left[K'_{2,\widehat\rho_{k_\rho}}(t)-K'_{2,\rho}(t)\right]dt \right|}$
\begin{eqnarray*}
&\leq&  \int_0^1 \left|t^{-1}W(t)-W(1)\right|t\left|K'_{2,\widehat\rho_{k_\rho}}(t)-K'_{2,\rho}(t)\right| dt\\
&\leq& |\widehat\rho_{k_\rho}| (1-\widehat\rho_{k_\rho}) \int_0^1  
\left|t^{-1}W(t)-W(1)\right| \left|t^{-\widehat\rho_{k_\rho}}-t^{-\rho}\right|dt
+|\widehat\rho_{k_\rho}-\rho|(1-\widehat\rho_{k_\rho}-\rho)\int_0^1 \left|t^{-1}W(t)-W(1)\right|dt\\
&=&o_\mathbb P(1).
\end{eqnarray*}
This achieves the proof of Theorem~2. $\blacksquare$

{\bf Proof of Corollary 4.} According to Theorem 2.1 in \cite{GdHP2002}, $\widehat \rho_{k_\rho}$ is consistent in probability as soon as the intermediate sequence $k_\rho$ satisfies $\sqrt {k_\rho} A\left({n \over k_\rho}\right) \to \infty$ and the second order condition $(C_{SO})$ hold. Combining this result with our Theorem~2, Corollary 4 follows. $\blacksquare$

{\bf Proof of Theorem~3.}  It is equivalent to show the asymptotic normality of
\begin{eqnarray*}
{\sqrt k \over \log {k\over np}} \log {\widehat x_{p, \xi} \over x_p}&=&
{\sqrt k \over \log {k\over np}} \left\{\log X_{n-k,n}+\widehat \gamma_k(K_{\widehat \Delta^*_{opt}}) \log {k \over np} -\log x_p \right.\\
&& \hspace{2cm}\left.- {(1-\xi)(1-2\xi)\over \xi^2} \left[\widehat \gamma_k(K_1)-\widehat \gamma_k(K_{2,\xi})\right] {\left({k \over np}\right)^{\xi}-1 \over \xi} \right\}\\
&=& \sqrt k \left(\widehat \gamma_k(K_{\widehat \Delta^*_{opt}})-\gamma\right) +{\sqrt k \over \log {k\over np}} \log {Q_n(1)\over U({n\over k})} - {\sqrt k \over \log {k\over np}} \left\{\log {U({1\over p})\over U({n\over k})} - \gamma\log {k\over np}\right\}\\
&& - {(1-\xi)(1-2\xi)\over \xi^2} {\sqrt k \left[\widehat \gamma_k(K_1)-\widehat \gamma_k(K_{2,\xi})\right] \over 
\log {k\over np}} {\left({k \over np}\right)^{\xi}-1 \over \xi} \\
&=&  \sqrt k \left(\widehat \gamma_k(K_{\widehat \Delta^*_{opt}})-\gamma\right) +{\sqrt k \over \log {k\over np}} \log {Q_n(1)\over U({n\over k})} - {\sqrt k \over \log {k\over np}} \widetilde A\left({n\over k}\right) {({k\over np})^\rho - 1\over \rho}\\
&&- {\sqrt k \over \log {k\over np}} \widetilde A\left({n\over k}\right) \left\{{\log U({1\over p})-\log U({n\over k})-\gamma \log {k\over np} \over \widetilde A({n\over k})} - {({k\over np})^\rho - 1\over \rho}\right\}\\
&& - {(1-\xi)(1-2\xi)\over \xi^2} {\sqrt k \left[\widehat \gamma_k(K_1)-\widehat \gamma_k(K_{2,\xi})\right] \over 
\log {k\over np}} {\left({k \over np}\right)^{\xi}-1 \over \xi} \\
&=:& T_3+T_4-T_5-T_6-T_7.
\end{eqnarray*}
For this aim, we will study the five terms separately. According to Theorem~2, we have
\begin{eqnarray*}
T_3 \stackrel{d}{\longrightarrow} {\cal N}(0, {\cal AV}(K_{\Delta^*_{opt}})).
\label{Q1}
\end{eqnarray*}
Now, according to Proposition 1, we have almost surely
\vskip2ex
\noindent
$\displaystyle{\left|\sqrt k \log {Q_n(1)\over U({n\over k})} - \gamma W(1)\right| }$
\begin{eqnarray*}
\leq \sup_{t\in (0, 1]} t^{{1\over 2}+\varepsilon} \left|\sqrt k \left(\log {Q_n(t)\over U({n \over k})} +\gamma \log t\right) - \gamma t^{-1} W(t) - \sqrt k \widetilde A\left({n \over k}\right) {t^{-\rho}-1\over \rho}\right| = o(1),  
\end{eqnarray*}
from which we deduce that
\begin{eqnarray*}
T_4 \stackrel{\mathbb P}{\longrightarrow} 0.
\label{Q2}
\end{eqnarray*}
Clearly, under our assumptions, we also have 
\begin{eqnarray*}
T_5 \longrightarrow 0.
\label{Q3}
\end{eqnarray*}
Now, according to the inequality (\ref{inequality})
\begin{eqnarray*}
\nonumber
|T_6|&\leq &{\sqrt k |\widetilde A({n \over k})| \over \log{k\over np}} \, \left|{\log U({1\over p})-\log U({n \over k})-\gamma\log {k\over np}\over \widetilde A({n \over k})} - {({k\over np})^\rho - 1\over \rho}\right|\\
\nonumber
 &\leq & {\sqrt k |\widetilde A({n \over k})| \over \log{k\over np}} \, \varepsilon \left({k \over np}\right)^{\rho+\delta}\\
 &=&o(1)
 \label{Q4}
\end{eqnarray*}
for any $0<\delta<-\rho$.\\
Finally, to treat $T_7$ two cases have to be considered: either $\xi$ is a canonical negative value $\widetilde \rho$ or an estimator $\widehat \rho$ consistent in probability such that $|\widehat \rho-\rho|=O_\mathbb P(n^{-\varepsilon})$ for some $\varepsilon>0.$\\
If a canonical negative value $\widetilde \rho$ is used, then according to Corollary 1, two times applied, we have
$$\sqrt k \left[\widehat \gamma_k(K_1)-\widehat \gamma_k(K_{2,\widetilde \rho})\right]=O_\mathbb P(1).$$ 
This immediately implies that $T_7=o_\mathbb P(1).$\\
If an estimator $\widehat \rho$ consistent in probability is used, we consider the decomposition
\begin{eqnarray*}
T_7&=& {(1-\widehat \rho)(1-2\widehat \rho)\over \widehat \rho^2} {\sqrt k \left[\widehat \gamma_k(K_1)-\widehat \gamma_k(K_{2,\widehat \rho})\right] \over 
\log {k\over np}} {\left({k \over np}\right)^{\rho}-1 \over  \rho}\\
&&+{(1-\widehat \rho)(1-2\widehat \rho)\over \widehat \rho^2} {\sqrt k \left[\widehat \gamma_k(K_1)-\widehat \gamma_k(K_{2,\widehat \rho})\right] \over 
\log {k\over np}} \left\{{\left({k \over np}\right)^{\widehat \rho}-1 \over  \widehat \rho}-{\left({k \over np}\right)^{\rho}-1 \over  \rho}\right\}.
\end{eqnarray*}
Note that
\begin{eqnarray*}
\sqrt k \left[\widehat \gamma_k(K_1)-\widehat \gamma_k(K_{2,\widehat \rho})\right]&=&\sqrt k \left[\widehat \gamma_k(K_1)-\gamma\right]-\sqrt k \left[\widehat \gamma_k(K_{2,\rho})-\gamma\right]-\sqrt k \left[\widehat \gamma_k(K_{2,\widehat \rho})-\widehat \gamma_k(K_{2,\rho})\right]\\
&=&O_\mathbb P(1),
\end{eqnarray*}
where we use Corollary 1 for the two first-terms of the right-hand side  and the proof of Theorem~2 (term $T_2$) for the last term. This implies that
\begin{eqnarray*}
T_7&=& o_\mathbb P(1)+o_\mathbb P(1) \left\{{\left({k \over np}\right)^{\widehat \rho}-1 \over  \widehat \rho}-{\left({k \over np}\right)^{\rho}-1 \over  \rho}\right\}\\
&=& o_\mathbb P(1)+o_\mathbb P(1) \int_1^{k/(np)} s^{\rho-1} \left[s^{\widehat \rho-\rho}-1\right] ds.
\end{eqnarray*}
Inspired by \cite{dHR1993}, we study this integral by using the inequality
$$\left|{e^x-1 \over x}-1\right| \leq e^{|x|}-1, \quad \forall x\in \mathbb R,$$
from which we deduce that
\vskip2ex
\noindent
$\displaystyle{\left|\int_1^{k/(np)} s^{\rho-1} \left[s^{\widehat \rho-\rho}-1\right] ds - \left(\widehat \rho-\rho\right) \int_1^{k/(np)} s^{\rho-1} \log s \, ds\right|}$
\begin{eqnarray*}
&\leq&\left| \left(\widehat \rho-\rho\right) \int_1^{k/(np)} s^{\rho-1} \log s \left\{{\exp\{\left(\widehat \rho-\rho\right)\log s\}-1 \over \left(\widehat \rho-\rho\right)\log s}-1 \right\} ds\right|\\
&\leq& \left| \widehat \rho-\rho\right| \int_1^{k/(np)} s^{\rho-1} \log s \left\{ \exp(\left|\widehat \rho-\rho\right|\log s)-1 \right\} ds\\
&\leq& \left| \widehat \rho-\rho\right|  \left\{\exp\left(\left|\widehat \rho-\rho\right|\log {k \over np}\right)-1  \right\}\int_1^{k/(np)} s^{\rho-1} \log s \, ds\\
&=&O_\mathbb P\left(\left( \widehat \rho-\rho\right)^2 \log{k \over np}  \int_1^{k/(np)} s^{\rho-1} \log s \, ds\right)\\
&=&o_\mathbb P(1).
\end{eqnarray*}
 This implies that
 \begin{eqnarray*}
 \int_1^{k/(np)} s^{\rho-1} \left[s^{\widehat \rho-\rho}-1\right] ds = \left(\widehat \rho-\rho\right) \int_1^{k/(np)} s^{\rho-1} \log s \, ds+o_\mathbb P(1)=o_\mathbb P(1).
\end{eqnarray*} 
 This entails $T_7=o_\mathbb P(1)$ and thus achieving the proof of Theorem~3. $\blacksquare$

\bibliography{mybib} 
\section{References}
B\"uhlmann, P. (2002), Bootstraps for time series, {\it Statistical Science}, {\bf 17(1)}, 52--72.\\
Cont, R. (2005), Volatility Clustering in Financial Markets: Empirical Facts and Agent–Based Models, {\it Long memory in economics}, A Kirman \& G Teyssiere (eds.): Springer.\\
Chavez-Demoulin, V., Embrechts, P., Sardy, S. (2014), Extreme-quantile tracking for financial time series, {\it Journal of Econometrics}, {\bf 181(1)}, 44--52.\\
Danielsson, J. (2011), {\it Financial Risk Forecasting}. Wiley.\\
Davison, A.C. and Hinkley, D.V. (1997), {\it Bootstrap Methods and Their Applications}, Cambridge Univ. Press. \\
Drees, H. (2000), Weighted approximations of tail processes for $\beta$-mixing random variables, {\it Ann. Appl.
Probab.}, {\bf 10}, 1274--1301.\\
Drees, H. (2003), Extreme quantile estimation for dependent data, with applications to finance, {\it Bernoulli}, {\bf 9},
617--657.\\
Embrechts, P., Koch, E., Robert, C. (2016), Space-time max-stable models with spectral separability, {\it  Probability, Analysis and Number Theory. In honour of N.H. Bingham}, C.M. Goldie and A. Mijatovi\c c (Eds.), {\it Advances of Applied Probability} Special Vol. {\bf 48A}, 77--97.\\
Gardes, L. and Girard, S. (2008), A moving window approach for nonparametric estimation of the conditional tail index, {\it J. Multivariate Anal.}, {\bf 99}, 2368-2388.\\
Goegebeur, Y. and Guillou, A. (2013), Asymptotically unbiased estimation of the coefficient of tail dependence, {\it Scand. J. Stat.}, {\bf 40}, 174--189.\\
Gomes, M.I., de Haan, L. and Peng, L. (2002), Semi-parametric estimation of the second order parameter in statistics of extremes, {\it Extremes}, {\bf 5}, 387--414.\\
de Haan, L. and Ferreira, A. (2006). {\it Extreme Value Theory. An Introduction}, Springer Series in Operations
Research and Financial Engineering, Springer, New York.\\
de Haan, L., Mercadier, C. and Zhou, C. (2016), Adapting extreme value statistics to financial time
series: dealing with bias and serial dependence, {\it Finance Stoch.}, {\bf 20}, 321--354.\\
de Haan, L. and Rootz\'en, H. (1993), On the estimation of high quantiles, {\it J. Statist. Plann. Inference}, {\bf 35}, 1--13.\\
Hill, B. (1975), A simple general approach to inference about the tail of a distribution, {\it Ann. Statist.}, {\bf 3}, 1163--1174.\\
Hsing, T. (1991), On tail index estimation using dependent data, {\it Ann. Statist.}, {\bf 19}, 1547--1569.\\
Kupiec, P.  (1995), Techniques for verifying the accuracy of risk management models. {\it J. of Derivatives}. {\bf 3}, 73–-84.\\
Liu, H., Shi, J. and Erdem, E. (2013), An integrated wind power forecasting methodology:
interval estimation of wind speed, operation probability of wind
turbine, and conditional expected wind power output of a wind farm, {\it International Journal of Green Energy}, {\bf 10}, 151--176.\\
Lojowska, A.,  Kurowicka, D., Papaefthymiou, G. and van der Sluis, L. (2010), Advantages of ARMA-GARCH wind speed time series modeling. {\it 2010 IEEE 11th International Conference on Probabilistic Methods Applied to Power Systems}, 83--88.\\
Matthys, G., Delafosse, E., Guillou, A. and Beirlant, J. (2004), Estimating catastrophic quantile levels for heavy-tailed distributions, {\it Insurance Math. Econom.}, {\bf 34}, 517--537.\\
McNeil, A.J. and Frey, R. (2000), Estimation of tail-related risk measures for heteroscedastic
financial time series: an extreme value approach, {\it Journal of Empirical Finance}, {\bf 7}, 271--300.\\
McNeil, A., Frey, R. and Embrechts, P. (2015), {\it Quantitative Risk Management:
Concepts, Techniques and Tools}. Revised Edition,  Princeton University Press.\\
St\u aric\u a, C. (1999), On the tail empirical process of solutions of stochastic difference equations, available at \url{http://143.248.27.21/mathnet/paper_file//Chalmers/Catalin/aarch.pdf}.\\
Tarullo, D.K. (2008), {\it Banking on Basel. The Future of International Financial Regulation}, Peterson Institute for International Economics, Washington, DC.\\
Swiss Re. (2014), Natural catastrophes and man-made disasters in 2013: large losses from floods and
hail; Haiyan hits the Philippines. Sigma 1/2014, Swiss Re, Zurich.\\
Weissman, I. (1978), Estimation of parameters and large quantiles based on the $k$ largest observations, {\it J. Amer. Statist. Assoc.}, {\bf 73}, 812--815.\\
\end{document}

%% file: colors.tex
\usepackage{color}



%% file: environments.tex
\newtheorem{theorem}{Theorem}


%% file: header.tex
\usepackage{url}
\usepackage{amsmath,amsfonts,amssymb,amscd,mathtools,bm,bbm}
\usepackage{multirow,tabularx,subfig,graphicx}
\usepackage{tikz}
\usetikzlibrary{calc,trees,positioning,arrows,chains,shapes.geometric, decorations.pathreplacing,decorations.pathmorphing,shapes, matrix,shapes.symbols}
\tikzset{ >=stealth', punktchain/.style={ rectangle, rounded corners, draw=black, very thick, text width=10em, minimum height=3em, text centered, on chain}}


%% file: tools.tex
\usepackage{bm}

\usepackage{amsmath}

\usepackage{mathrsfs}

